\documentclass[12pt]{article}
\usepackage{amsfonts}

\title{ Existence results for parabolic problems related to fully non linear operators degenerate or singular }

\author{ F. Demengel\\
University of Cergy-Pontoise, site de saint martin, 95302, Cergy-Pontoise\\
email  : demengel@math.u-cergy.fr\\
fax number : (33) 134256546}

\date{}

\newtheorem{theo}{Theorem}
\newtheorem{prop}{Proposition}
\newtheorem{rema}{Remark}
\newtheorem{defi}{Definition}
\newtheorem{cor}{Corollary}
\newtheorem{lemme}{Lemma}

\def\R{{\rm I}\!{\rm  R}}

\newcommand{\N}{{\bf N}}

\begin{document}

 \maketitle 
 \begin{abstract}

 In this paper we prove some existence and regularity results concerning parabolic equations 
 $$u_t = F(\nabla  u, D^2 u) + f(x,u)$$
  with some boundary conditions , on $\Omega \times ]0,T[$, where $\Omega$ is some bounded  domain which possesses the cone property and $F$ is singular or degenerate, with some uniform ellipticity conditions.
  
  \end{abstract}
\section{Introduction and hypothesis}

In this paper we  consider  the  parabolic equation 

$$
u_t= F(x, \nabla u, D^2 u) + h(x,t)\cdot \nabla u|\nabla u|^\alpha + f(x,t)$$
on some bounded domain $ Q_T = \Omega\times ]0,T[$ of $\R^N$, with some non zero  boundary conditions on the parabolic boundary.

Here the  operator is fully non linear and degenerate or singular, it satisfies some  assumptions as  in \cite{BD1}, which will be detailed later. In particular the class of operators considered contains both the $p$-Laplace and the Pucci operators, as well as non variational extensions of the $p$-Laplacian. Both $h$ and $f$ are bounded and continuous functions.  The boundary condition will be supposed to be  H\"older continuous. 

In previous papers \cite{BD1,BD2,BD5},the author, in collaboration with Birindelli, has considered the stationary case, introducing the notion of principal eigenvalue and proving the existence of solutions for a large class of Dirichlet problems. The parabolic case treated here requires the introduction of many new tools and new ideas.

 We begin by stating a definition of viscosity solutions adapted to the context, the difficulty being that due to the fact that the operator $F$ is not defined  when the gradient is zero, one cannot test points on which every test function have the gradient equal to zero. In the stationary case this is solved by just "not testing" such points unless the solution is locally constant. Here the situation is more involved and requires some "testing".

 The key points to prove the existence of solution are on one hand, some comparison principle and on the second hand the existence of some upper and lower barriers. 
 
 The comparison principle   presents some difficulty linked to the non definition of the operators when the gradient of test functions is zero, difficulty overcome  with the aid of the adapted definition of viscosity solutions that we propose.  This  comparison theorem  permits in particular to get the uniqueness of solution.
 
 The existence of lower and upper barriers is complicated by the  fact that the operator is homogeneous with different powers  with respect to $t$ and $x$,  a difference with most of  the  papers cited before.  
 
   In a third time, we use  Perron's method adapted to the context.

  We also establish some  regularity result, more precisely the solutions are H\"older in both the spatial and the time variable,  with some exponent which depends on the regularity of the  data $f$ and of the boundary value $\psi$, and also on the parameters of the exterior cone related to the open set $\Omega$. 
  
 Finally we also consider  to the case  of  some infinite domain such as 
   $\Omega \times \R^+$ and $\R^N\times ]0,T[$. 
     
\bigskip
Analogous problems are studied by Crandall, Kocan, Lions , and Swiech in  \cite{CKLS}  for  the case of Pucci's operators,  by Ishii and Souganidis  \cite{IS} for operators singular or degenerate and homogeneous of degree 1, and by Onhuma and  Sato  \cite{OS}  in the case of the $p$-Laplacian.

In \cite{JK}  and \cite{ES} , Juutinen and Kawhol treat the case of the infinite  Laplacian
 when the right hand side  $f $ is zero  and  the open domain is regular. Let us note that this situation is analogous to the present  one when   $\alpha = 0 $.  In their situation the  operator is linear   with respect to $D^2 u$ but  it is not well defined on points where the gradient  is zero. This   leads the authors to  give a convenient definition of viscosity solutions. This definition provides  a comparison principle and in particular the solutions  obtained are unique . The existence is obtained through a regularizing process, and using classical  results of Ladishenskaia Uralceva for parabolic problems. 
  
  On the other hand  \cite {CKLS} the authors consider the case of Pucci's operators in domains which have the uniform exterior cone condition, and with a right hand side $f$ bounded. They exhibit a supersolution and a sub-solution constructed with the aid of the  parameters of the cone relative to $\Omega$. They also prove a comparison principle which enables them to prove that the sub-solution is less than the supersolution. Finally through the Perron's method they prove the existence of a solution. 
  
  In \cite {OS} the authors consider the case of the $p$-Laplacian and  a right hand side   zero. They give a convenient definition of viscosity solution which provides a comparison principle . 
   This definition of viscosity solutions requires to introduce a set of admissible test functions when the gradient of $u$ is zero. Since  it  can be extended to our situation, it is natural to check that it is  equivalent to our definition, which is  done   in the appendix.

 \section{Notations and hypothesis}
In all that paper, (except in section 6) we shall assume that   $\Omega$ is some bounded domain  which   satisfies the  uniform exterior cone condition, .e. we assume that  there exist $\phi\in ]0, \pi[$ and $\bar r>0$ such that    for any $z\in\partial \Omega$ and  for an axe through $z$ of direction $\vec n_z$,
 $$ T_\phi=\{x:\ \frac{(x-z)\cdot\vec n_z}{|z-x|}\leq \cos\phi\}, \quad\quad T_\phi \cap\overline{\Omega}\cap B_{\bar r}(z)=\{z\}.$$

For a real $T$ positive  let $Q_T = \Omega \times ]0,T[$. We shall denote by $\partial Q_T $ the parabolic  boundary $(\partial \Omega \times ]0,T[) \cup (\Omega\times \{0\})$.
 Concerning $F$ we shall assume that $\alpha >-1$ and $F$ satisfies

\begin{itemize}

\item[(H1)] 
 $F: \Omega\times \R^N\setminus\{0\}\times S\rightarrow\R$,  is continuous with respect to  all its variables, 
and  $\forall t\in \R^\star$, $\mu\geq 0$, for all $x\in \Omega$, $p\neq 0$ and $X\in S$,  
 $$F(x, tp,\mu X)=|t|^{\alpha}\mu F(x, p,X).$$

\item[(H2)]
For $x\in \overline{\Omega}$, $p\in \R^N\backslash \{0\}$, $M\in S$,  $N\in S$, 
$N\geq 0$

\begin{equation}\label{eqaA}
a|p|^\alpha tr(N)\leq F(x,p,M+N)-F(x,p,M) \leq A
|p|^\alpha tr(N).
\end{equation}

\item [(H3)]
 There exists a
continuous function $  \omega$ with $\omega (0) = 0$, such that if
$(X,Y)\in S^2$ and 
$\zeta\in \R^+$ satisfy
$$-\zeta \left(\begin{array}{cc} I&0\\
0&I
\end{array}
\right)\leq \left(\begin{array}{cc}
X&0\\
0&Y
\end{array}\right)\leq 4\zeta \left( \begin{array}{cc}
I&-I\\
-I&I\end{array}\right)$$
and $I$ is the identity matrix in $\R^N$,
then for all  $(x,y)\in \R^N$, $x\neq y$
$$F(x, \zeta(x-y), X)-F(y,  \zeta(x-y), -Y)\leq \omega
(\zeta|x-y|^2).$$

  Sometimes this condition $(H3)$  can be replaced by the weaker assumption, which will for example be employed to  prove Holder's regularity results : 
 \item[(H4)]
There exists a continuous function $\tilde \omega$, $\tilde \omega(0)=0$
such that for all $x,y,$ in $\Omega$,  $p\neq 0$, $\forall X\in S$
$$|F(x,p,X)-F(y,p,X)|\leq \tilde \omega(|x-y|) |p|^\alpha |X|.$$
\end{itemize}

We assume that $h$ is continuous and bounded on $Q_T$  with  values in $\R^N$  and  satisfies (H5) :

There exists $\omega_h\leq 1$ and $c_h>0$  such that for all $(x,t)$, $(x,s)$ in $Q_T$
$$|h(x,t)-h(x,s)|\leq c_h|t-s|^{\omega_h }.$$

 Furthermore 

- Either $\alpha \leq 0$  
and  for all $(x,y)$ in $\Omega$ and $t\in ]0,T[$ 
$$|h(x,t)-h(y,t)|\leq c_h|x-y|^{1+\alpha}$$

- or  $\alpha >0$ and  for all $(x,y)$ in $\Omega$ and $t\in ]0,T[$ 

$(h(x,t)-h(y,t)\cdot x-y)\leq 0$. 

Concerning $f$ we shall assume that it  is at least continuous and will precise  further regularity when it will be needed.

  \bigskip
We now give the definition of viscosity solutions  adapted to our context. 

 It is well known that when  dealing with viscosity respectively  sub and super
solutions one works with  

$$u^\star (x,t) = \limsup_{(y,\tau), |(y,\tau)-(x, t)|\leq r} u(y,\tau)$$
and 
$$u_\star (x,t) = \liminf_{(y,\tau), |(y,\tau)-(x, t)|< r}\ u(y,\tau).$$
 It is easy to see that 
$u_\star \leq u\leq u^\star$ and $u^\star$ is upper semicontinuous (USC),
$u_\star$ is lower semicontinuous (LSC).  See e.g. \cite{CIL, I}.

 \begin{defi}

We shall say that $u$,   locally bounded, is a viscosity subsolution of 
$$u_t-F(x, \nabla u , D^2 u)- h(x,t) \cdot \nabla u|\nabla u|^\alpha \leq f(x,t)\quad \mbox{in}\quad \Omega\times (0,T)$$
if , for any $(\bar x,\bar t)\in \Omega\times (0,T)$, 

\begin{itemize}
\item either for all $\varphi\in  {\cal C}^2$ touching  $u^\star$ by above  at $\bar x$ such that 
  $\nabla_x \varphi(\bar x, \bar t)\neq  0$
$$\varphi_t(\bar x, \bar t) - F(\bar x,\nabla \varphi(\bar x, \bar t),  D^2\varphi(\bar x, \bar t))- h(\bar x, \bar t) \cdot \nabla \varphi|\nabla \varphi|^\alpha(\bar x, \bar t)\leq f(\bar x, \bar t).$$

\item or, if there exists $\delta_1$ and   $\varphi\in{\cal C}^2(]\bar t-\delta_1,\bar t+\delta_1[)$,  such that for any $t\in  ]\bar t-\delta_1,\bar t+\delta_1[$ 
$$\left\{\begin{array}{l}
\varphi(\bar t)=0\\
 u^\star (\bar x, \bar t)\geq  u^\star (\bar x,t)-\varphi(t)\\
\displaystyle \sup_{t\in ]\bar t-\delta_1,\bar t+\delta_1[} (u^\star(x,t)-\varphi(t)) \mbox{ is constant in a neighborhood of  }\ \bar x,
\end{array}
\right.
$$
  then 
 $$\varphi^\prime ( \bar t)\leq f(\bar x, \bar t). $$
 
 \end{itemize}

$u$,  locally bounded, is a viscosity supersolution of 
$$u_t -F(x,\nabla u, D^2 u)- h( x, t) \cdot \nabla u|\nabla u|^\alpha\geq f\quad \mbox{in}\quad \Omega\times (0,T)$$
if , for any $(\bar x,\bar t)\in \Omega\times (0,T)$, 

\begin{itemize}
\item either  for all $\varphi\in  {\cal C}^2$ which touches $u_\star $ by below   at $\bar x$,  
 such that  $\nabla_x \varphi(\bar x, \bar t)\neq  0$, 
$$\varphi_t (\bar x, \bar t)-F(\bar x, \nabla \varphi(\bar x, \bar t),  D^2\varphi(\bar x, \bar t))- h(\bar x, \bar t) \cdot \nabla \varphi|\nabla \varphi|^\alpha (\bar x, \bar t)\geq f(\bar x, \bar t).$$

\item or, if  there exists $\delta_1$ and   $\varphi\in{\cal C}^2(]\bar t-\delta_1,\bar t+\delta_1[)$ 
such that for any $t\in  ]\bar t-\delta_1,\bar t+\delta_1[$ 
$$\left\{\begin{array}{l}
\varphi(\bar t)=0\\
 u_\star(\bar x, \bar t)\leq  u_\star (\bar x,t)-\varphi(t)\\
\displaystyle \inf_{t\in ]\bar t-\delta_1,\bar t+\delta_1[} (u_\star(x,t)-\varphi(t)) \mbox{ is locally constant in a neighborhood of  }\ \bar x,
\end{array}
\right.
$$
then
 $$\varphi^\prime ( \bar t)\leq f(\bar x, \bar t). $$

\end{itemize}

Finally  a continuous function $u$ is a viscosity solution when $u$ is both a viscosity sub and supersolution. 
\end{defi}
 \begin{rema}
  In the following and for convenience of the reader we recall the definition of semi-jets  for parabolic problems : 
  $$ J^{2,+}u(\bar x, \bar t ) = \{ ( q, p, X)\in \R \times \R^N\times S,\ q(t-\bar t) + p. (x-\bar x) +{1\over 2} ^t(x-\bar x)X (x-\bar x) \geq u(x, t)-u(\bar x, \bar t)\}$$
and 
  $$ J^{2,-}u(\bar x, \bar t ) = \{ ( q, p, X)\in \R \times \R^N\times S,\ q(t-\bar t) + p. (x-\bar x) +{1\over 2} ^t(x-\bar x)X (x-\bar x) \leq u(x, t)-u(\bar x, \bar t)\}$$
\end{rema}

\begin {rema}
We prove in the appendix that our solutions are the same as those of Onhuma and Sato in the case where $\alpha \neq 0$, and to those of Evans and Spruck and Juutinen and Kawohl in  the case of the infinity Laplacian.  

\end{rema}

In the following we shall denote by $1_{\{f\}}$ the equation 
$$
u_t= F(x, \nabla u, D^2 u) + h(x,t)\cdot \nabla u|\nabla u|^\alpha + f(x,t)$$
 and by $1_{\{f, \psi\}}$ the boundary value problem 
 $$\left\{ \begin{array}{cc}
u_t= F(x, \nabla u, D^2 u) + h(x,t)\cdot \nabla u|\nabla u|^\alpha + f(x,t)&{\rm in}\ Q_T\\
 u(x,0) = \psi (x) &{\rm on}\  \partial Q_T 
 \end{array}\right.$$

\begin{rema}

 Let us note that if $u$ is a sub-solution (respectively supersolution) of $1_{f}$ and if $\varphi$ is some  ${\cal C}^1$  function  depending only on $t$,  $(x,t)\mapsto u(x,t)+ \varphi (t)$ is a sub-solution (respectively supersolution) of 
 $1_{\{f+ \varphi^\prime \}}$. 
 \end{rema}

\section{Comparison principle and barriers.}
We begin to prove a  comparison principle for the operator $u_t-F(x, \nabla u, D^2 u)-h(x)\cdot \nabla u|\nabla u|^\alpha$.  One of its consequences is  the uniqueness of the solutions for $1_{\{f, \psi\}}$. 
 
 \begin{theo}\label{thcomp} Suppose that $u$  is a sub-solution bounded  for $1_{\{g\}}$ and  $v$  is   a  super solution  bounded  of  $1_{\{ f\}}$ with $g\leq  f$ in $\Omega\times ]0,T[$,  $g$  being upper semicontinuous  and $f$ being   lower semicontinuous. Suppose that $u^\star\leq v_\star$  on  $(\partial\Omega\times [0,T))\cup(\Omega\times\{0\})$,  
then $u^\star\leq v_\star $ in $\Omega\times [0,T)$.
\end{theo}

The proof of this theorem requires the following technical lemma which proof is postponed after the proof of theorem \ref{thcomp}  for the sake of clearness. 

 \begin{lemme}\label{lem1} 
 Suppose that $\Omega$ is some  open set. 
 Suppose that $u$ is a supersolution of 
 $$u_t-F(x, \nabla u, D^2 u)-h(x,t)\cdot \nabla u|\nabla u|^\alpha \geq f(x,t)$$
in $Q_T = \Omega \times ]0,T[$  and suppose that $C_1$   is  some  constant , that $\varphi$ is some ${\cal C}^2$ function on $]0,T[$,  that $k > \sup (2, {\alpha+2\over \alpha+1})$
 and  $(0, \bar t)\in \Omega \times ]0, T[$ are  such that  for some  $\delta_1>0$
 $$\inf _{x\in B(0, \delta_1), | t-\bar t|< \delta_1}  (u(x,t)-\varphi(t)+C_1 |x|^k) =  u(0, \bar t)$$
 Then 
 $$\varphi^\prime (\bar t) \geq f(0, \bar t).$$
  
 \end{lemme}

Proof of theorem \ref{thcomp}   : 
          
               Suppose by contradiction that 
               $u(\bar x, \bar t) > v(\bar x, \bar t)$ for some $(\bar x, \bar t)\in Q_T$, let $\kappa>0$ be such that 
               $${2\kappa\over T-\bar t} < {(u-v)(\bar x, \bar t)\over 2}, $$
                then 
                $u_1(x,t)= u(x,t)-{\kappa\over T-t} $ is a strict sub-solution, $v_1(x,t) = v(x,t)+ {\kappa\over T-t}$ is a strict supersolution and $u_1-v_1>0$ somewhere in $Q_T$. Moreover the maximum of $u_1-v_1$  cannot be achieved  in $|t-T| < {T-\bar t\over 2}$, since in that set one has 
                $$u-{K\over T-t}-(v+ {K\over T-t} )\leq \sup (u-v)-{4\kappa \over T-\bar t}$$
                while 
                $$u(\bar x, \bar t)-{K\over T-\bar t}-(v(\bar x, \bar t)+ {K\over T-\bar t})\geq   \sup (u-v) -{2\kappa\over T-\bar t}. $$
In the following we replace $u$ by $u-{\kappa\over T-t}$ which is a   sub-solution of 
$1_{f-{\kappa\over (T-t)^2}}$ and $v$  by $v+ {\kappa\over T-t}$ which a supersolution of $1_{f+ {\kappa\over (T-t)^2}}$. 
                
                 We  define for $j\in \N$  and for $k> \sup (2, {\alpha+2\over \alpha+1}, {2(1+\alpha)\over \omega_h})$, 
                 $$\Psi_j(x,t,y,s) = u^\star(x,t)-v_\star (y,s)-{j\over 2} |t-s|^2-{j\over k} |x-y|^k$$
Then  $\psi_j$ achieves its maximum on $(x_j, t_j, y_j, s_j)\in (\Omega\times ]0,T[)^2$. It is classical that the sequences
  $(x_j,t_j)$ $(y_j,s_j)$   both converge to $(\bar x,\bar t)$ which is a maximum point for $u^\star -v_\star$,  and that  $j|s_j-t_j|^2 + j|x_j-y_j|^k\rightarrow 0$. 
 
 We want to prove that for $j$ large enough $x_j\neq y_j$. Suppose not  i.e.  $x_j = y_j$ 
 then 
 $$(y,s)\mapsto  v_\star(x_j, s_j)-{j\over k} |x_j-y|^k-{j\over 2} |s-t_j|^2 + {j\over 2} |t_j-s_j|^2$$
 would be  a test function from below for $v_\star $ at $(x_j, s_j)$.  Then  
   applying Lemma \ref{lem1}   in its form for super-solutions with $C_1 = {j\over k}$,  $\varphi $  replaced by $t\mapsto  v_\star(x_j, s_j)-{j\over 2} |t-t_j|^2 + {j\over 2} |t_j-s_j|^2$,   replacing $0$ by $x_j$,  and $\bar t$ by $s_j$ one would get  that

  $$-j(s_j-t_j)\geq  {\kappa\over T^2}+f(x_j, s_j).$$
  On the other hand 
  $$(x,t)\mapsto u^\star (x_j, t_j) + {j\over k} |x_j-x|^k + {j\over 2} |t-s_j|^2-{j\over 2} |t_j-s_j|^2$$
   would be a test function from above for $u^\star$ on $(x_j, t_j)$.  
 Using Lemma \ref{lem1} in its form  for sub-solutions, with $\varphi
 $ replaced by  
 $t\mapsto u(x_j, t_j) + {j\over 2} |t-t_j|^2-{j\over 2} |t_j-s_j|^2$ 
  $0$ by $x_j$,   $C_1$ by $-{j\over k}$, one gets that 
 $$j(t_j-s_j) \leq g(x_j,t_j)-{\kappa\over T^2}.$$
Substracting the two inequalities,  passing to the limit  and using the upper semicontinuity of $g$ and the  lower semicontinuity of $f$,  one gets  that
 $$\lim_{j\rightarrow+\infty}j(t_j-s_j)+ j(s_j-t_j)\leq -{2\kappa\over T^2}+ \limsup_{j\rightarrow +\infty} (g(x_j, t_j)-f(x_j, s_j))\leq  -{2\kappa\over T^2}$$
which is  a contradiction. 
 
 We have then proved that $x_j\neq y_j$.

By Ishii's lemma, (see also lemma 2.1 in \cite{BD1})  there exist   
$(X_j,Y_j)\in S^2$ , with 
 $$(j(t_j-s_j), j|x_j-y_j|^{k-2}(x_j-y_j), X_j)\in J^{2,+} u^\star (x_j,  t_j)$$
  
  $$(j(t_j-s_ j), j|x_j-y_j|^{k-2}(x_j-y_j), -Y_j)\in J^{2,-} v_\star(y_j,s_j)$$
and for some  positive constant $c$ 
$$\left(\begin{array}{cc} X_j&0\\
0&Y_j\end{array}\right) \leq cj|x_j-y_j|^{k-2} \left(\begin{array}{cc} I&-I\\
-I&I
\end{array}\right)$$
This implies that, using assumption (H3) and the fact that $j|x_j-y_j|^k\rightarrow 0$
 \begin{eqnarray*}
{\kappa\over T^2}+ f(y_j, s_j)&\leq & j(t_j-s_j) -F(y_j, j|x_j-y_j|^{k-2}(x_j-y_j), -Y_j)\\
&+& j^{1+\alpha}h(y_j, s_j)
\cdot (x_j-y_j)|x_j-y_j|^{k-2+ (k-1)\alpha}\\
&\leq  &j(t_j-s_j) -F(x_j, j|x_j-y_j|^{k-2}(x_j-y_j),  X_j)+ o(1) \\
&+& j^{1+\alpha} h(x_j, t_j)\cdot (x_j-y_j)|x_j-y_j|^{k-2+ (k-1)\alpha} + o(1)\\
&\leq &  g(x_j, t_j) -{\kappa\over T^2} +o(1),
\end{eqnarray*}
Using the lower semicontinuity of $f$, the uppersemicontinuity of $g$  and letting $j\rightarrow +\infty$ we get a contradition.

In the previous inequalities we have used 
\begin{eqnarray*}
 |h(x_j, t_j)-h(x_j, s_j)|&&|x_j-y_j|^{(k-1)(1+\alpha)}j^{1+\alpha}\\
 & \leq&  c_h |t_j-s_j|^{\omega_h}j^{1+\alpha} |x_j-y_j|^{(k-1)(1+\alpha)}\\
 &\leq& (j|t_j-s_j|^2)^{\omega_h\over 2} (j|x_j-y_j|^k)^{(1+\alpha)(k-1)\over k} j^{{1+\alpha\over k}-{\omega_h\over 2}}\\
 &= &o(1)
 \end{eqnarray*}
 and when $\alpha <0$
 $$|h(x_j, s_j)- h(y_j, s_j)||x_j-y_j|^{(k-1)(1+\alpha)}j^{1+\alpha}\leq j^{1+\alpha} |x_j-y_j|^{k(1+\alpha)}= o(1).$$

{\em Proof of Lemma  \ref{lem1} }
First replacing if necessary $\varphi$ by $\varphi (t)+ C_2 |t-\bar t|^2$ for some constant $C_2>0$ and $C_1$ by some constant $> C_1$ one can assume that the infimum is strict in $x$ and $t$ separately. 

Clearly $\psi(x,t)=\varphi(t) -C_1|x|^k-C_2(t-\bar t)^2$ is a test function for $u$ in $(0,\bar t)$ but its  gradient with respect to $x$  is zero.  So we are going to prove that either the function $t\mapsto \varphi (t) + C_2 |t-\bar t|^2$ is a test function as in the second case of the definition of viscosity supersolution and then the conclusion of the Lemma is immediate. Or, if this is not the case, then it is possible to construct a sequence of points tending to $(0,\bar t)$ for which there exists a test function which  gradient   with respect to $x$  is different from zero,  but tend to zero. Then passing to the limit we get the required inequality.

Hence we suppose first that  the function $t\mapsto \varphi (t) - C_2 |t-\bar t|^2$ is as in the definition of viscosity supersolution i.e. we suppose that there exists  $\delta_1>0$, and   $\bar \delta>0$ such that   for all $x\in B(0, \bar \delta)$,   
 $$\inf_{|t-\bar t|< \delta_1}  \{v(x,t)-\varphi(t)+ C_2(t-\bar t)^2\}= \inf_{|t-\bar t|< \delta_1} \{v(0,t)-\varphi(t) + C_2(t-\bar t)^2\}.$$ 
  We claim that this  infimum is achieved on  $(0, \bar t).$ Indeed, the infimum is less or equal to $v(0, \bar t)$ and on the other hand it is  more than 
  $ \inf_{x\in B(0, \delta_1), |t-\bar t|< \delta_1} \{ v(x,t)+ C_1|x|^k-\varphi(t)+ C_2(t-\bar t)^2\}$ which equals $v(0, \bar t)$.
  
  Then the conclusion  given in that case in the definition of viscosity supersolution is that 
  $ \varphi^\prime (\bar t)  \leq f(0, \bar t)$. 
  
  \bigskip
  
   We now suppose that we are not in this situation i.e. that $x\mapsto  \inf_{|t-\bar t|< \delta_1}  v(x,t)-\varphi(t) + C_2|t-\bar t|^2$ is not constant  in a neighborhood of $\bar x$.  
   
       Recall that since the infimum is strict in $x$ and $t$ separately, for all $\delta >0$, $\delta < \delta_1$  there exists $\epsilon (\delta)>0$ such that 
    \begin{eqnarray*}
    \inf \left(\inf_{|t-\bar t|> \delta, x\in B(0, \delta_1)}\{ v(x,t)+\right.&& C_1|x|^k-\varphi(t)+ C_2(t-\bar t)^2\}, \\
    &&\left. 
    \inf_{|t-\bar t|> \delta_1, |x|> \delta} \{ v(x,t)+ C|x|^k-\varphi(t)+ C_2(t-\bar t)^2\}\right) \\
    &\geq &v(0, \bar t) + \epsilon (\delta).
    \end{eqnarray*}
    
  We now choose $\delta_2 \leq \inf ({\epsilon (\delta) \over 4 C_1 k (2\delta_1 )^{k-1} }, \delta) $.
      Then,  with that choice,  for all  $x\in B(0, \delta_2)$
      $$\inf_{y\in B(0, \delta_1),|t-\bar t|\leq  \delta_1} \{ v(y,t)+ C_1|x-y|^k-\varphi(t)+ C_2(t-\bar t)^2\}
      \leq v(0, \bar t) +{ \epsilon (\delta)\over 4}$$
      while 
         $$\inf_{|y|> \delta,|t-\bar t|\leq  \delta_1} \{ v(y,t)-\varphi(t)+  C_1|y-x|^k-b(t-\bar t) + C_2(t-\bar t)^2\}
      \geq v(0, \bar t) +{3 \epsilon (\delta)\over 4}.$$
      
      Moreover one also has 
        \begin{eqnarray*}
        &&\inf_{y\in B(0, \delta_1),|t-\bar t|> \delta} \{ v(y,t)-\varphi(t)+ C|x-y|^k+ C_2(t-\bar t)^2\}\\
    &  \geq& \inf_{y\in B(0, \delta_1),|t-\bar t|> \delta } \{ v(y,t) -\varphi(t) + C_1|y|^k+ C_2(t-\bar t)^2\}-{\epsilon (\delta)\over 4}\\
    & \geq& v(0, \bar t)+{3\epsilon (\delta)\over 4}.
    \end{eqnarray*} 
    This implies  that  for all $x\in B(0, \delta_2)$

      \begin{eqnarray}\label{eqeq001}
       && \inf_{y\in B(0, \delta_1),|t-\bar t|< \delta_1}\{ v(y,t)+ C_1|y-x|^k-\varphi(t)+ C_2(t-\bar t)^2\}\nonumber\\
   &&=   \inf_{y\in B(0, \delta), |t-\bar t |\leq \delta } \{ v(y,t)+ C_1|y-x|^k-\varphi(t) + C_2(t-\bar t)^2\}. 
       \end{eqnarray}
       
        Since  $x\mapsto  \inf_{|t-\bar t|< \delta_1}  \{v(x,t)-\varphi(t)+ C_2|t-\bar t|^2\}$ is not constant  in a neighborhood of $\bar x$, 
   there exist $(x_\delta, y_\delta ) \in B(0, \delta_2 )$ 
   $$ \inf_{|t-\bar t|< \delta_1}  \{v(x_\delta ,t)-\varphi(t)+  C_2|t-\bar t|^2\}>  \inf_{|t-\bar t|< \delta_1} \{ v(y_\delta ,t)-\varphi(t)+C_1 |x_\delta -y_\delta |^k+  C_2|t-\bar t|^2\}$$
Hence
  $$  \inf_{y\in B(0, \delta_1), |t-\bar t|< \delta_1} \{v(y,t)-\varphi(t)+C_1 |x_\delta -y|^k+  C_2|t-\bar t|^2\}$$ is achieved on some point $(z_\delta , t_\delta )$ with $z_\delta  \neq x_\delta $.  Indeed if it was achieved on $(x_\delta, t_\delta)$ for some $t_\delta$  one would have 
    \begin{eqnarray*}
      v(x_\delta ,t_\delta )-\varphi(t_\delta ) &+&  C_2|t_\delta -\bar t|^2\\
      &= & \inf_{y\in B(0, \delta_1),|t-\bar t|< \delta_1}  \{v(y,t)-\varphi(t)+C_1 |x_\delta-y|^k+  C_2|t-\bar t|^2\}\\
      &\leq&   \inf_{|t-\bar t|< \delta_1} \{ v(y_\delta ,t)-\varphi(t)+C_1 |y_\delta-x_\delta  |^k+  C_2|t-\bar t|^2 \}\\
      &<  &\inf_{|t-\bar t|< \delta_1} \{ v(x_\delta ,t)-\varphi(t)+ C_2|t-\bar t|^2\}\\
      &\leq& v(x_\delta ,t_\delta )-\varphi(t_\delta )+  C_2|t_\delta -\bar t|^2, 
      \end{eqnarray*}
      a contradiction. 
  Moreover using  (\ref{eqeq001}),   the infimum is achieved in $B(0, \delta) \times ]\bar t-\delta, \bar t+\delta[$.

 All this imply that   $(y,t)\mapsto v(z_\delta, t_\delta) + \varphi(t) -\varphi(\bar t_\delta) + C_1 |x_\delta -z_\delta|^k- C_1|x_\delta -y|^k+ C_2(t_\delta-\bar t)^2 -C_2 |t-\bar t|^2$ is  a test function for $v$ on $(z_\delta, t_\delta)$  and  since $v$ is a supersolution 
        \begin{eqnarray*}
       \varphi^\prime (t_\delta ) - 2C_2 (t_\delta-\bar t) &-&F(- C_1k|x_\delta -z_\delta |^{k-2} (z_\delta -x_\delta), X_\delta) \\
        &+ & k^{1+\alpha} |x_\delta -z_\delta |^{(k-1)(\alpha +1)-1}  h(z_\delta, t_\delta)\cdot  (z_\delta -x_\delta)\\
        &\geq &f(z_\delta, t_\delta)
        \end{eqnarray*}
  where $X_\delta=-D^2 (C_1|x_\delta -y|^k)\mid_{y=z_\delta }$.     
 We have finally obtained that
 $$\varphi^\prime (t_\delta) - 2C_2(t_\delta-\bar t)+  C_1^{1+\alpha}   |x_\delta -z_\delta |^{k(\alpha+1)-\alpha-2}  + |h|_\infty k ^{1+\alpha} (2\delta) ^{(k-1)(\alpha+1)}  \geq  f(z_\delta, t_\delta).$$
Using $x_\delta \in B(0, \delta_2)\subset B(0, \delta)$,  $z_\delta\in B(0, \delta)$ and $k> {\alpha+2\over \alpha+1}$, 
 $$\varphi^\prime (t_\delta) + o(1) \geq f(z_\delta, t_\delta).$$
        
 Letting $\delta$ go to zero,  and using the lower semicontinuity of $f$ one gets the result. This ends the proof of lemma \ref{lem1}.

\bigskip
 
We now construct a supersolution and a subsolution  for  $1_{\{f,\psi\}}$
  We recall that   in \cite{BD5} we constructed a global barrier for the stationary case:
\begin{prop}\label{propbarglob}
For all $z\in \partial \Omega$, there exists  some     function $W_z$  continuous on $\overline{\Omega}$, such that 
$W_z(z)=0$,  $W_z>0$ in $\Omega \setminus \{z\}$, which satisfies 

$$ F(x,\nabla W_z, D^2 W_z)+ h(x,t)\cdot \nabla W_z |\nabla W_z|^\alpha\leq -1\quad\mbox{in}\quad \Omega.$$
Furthermore  $\nabla W_z\neq 0$ everywhere and there exist $\underline{c}
>0,$  $\overline{c}>0$  and $\gamma \in ]0,1[$ which depend on the parameters of the cone, such that  for all $z\in \partial \Omega$ and $x\in \Omega $ 
$$\underline{c}|z-x|^\gamma \leq W_z(x) \leq \overline{ c}|x-z|^\gamma . $$

\end{prop}

\begin{rema} In fact one can ask, up to change the constants $\gamma$  and the constants $\underline{c}$ and $\overline{c}$  that  $W_z$  be such that $-W_z$  be also  a sub-solution of 
$$F (x,\nabla (-W_z) , D^2(- W_z) )- h(x,t)\cdot \nabla W_z  |\nabla W_z |^\alpha\geq 1\quad\mbox{in}\quad \Omega.$$

\end{rema}
The proof of Proposition \ref{propbarglob} can be found in \cite{BD5}.

\bigskip

We now give some existence's result of supersolutions and sub-solutions  for the parabolic problem. 

\begin{prop}\label{sub}
 Suppose that $\psi$ is Lipschitzian in $t$, H\"olderian with exponent 
 $\gamma$ in $x$. Suppose that $f$ is uniformly bounded. Then there exists a  continuous supersolution $W$ of 
$1_{\{|f|_\infty, \psi\}}$. 

In the same manner  there exists a continuous sub-solution $V$ of $1_{\{-|f|_\infty, \psi\}}$. 
\end{prop}
Proof of proposition \ref{sub}.

Let $c_\psi$ be some holder's constant for $\psi$.  We define

$$
  W_1(x,t) := \inf_{(z, \tau)\in \partial \Omega \times ]0,T[} \{ \psi (z, \tau)+ \left({c_\psi\over \underline{c}} 
  +( |\psi_t|_\infty + |f|_\infty )^{1\over 1+\alpha}\right) W_z(x) + |\psi_t|_\infty |t-\tau|\}.$$
  
  Let us note that 
  $ \left({c_\psi\over \underline{ c}}+( |\psi_t|_\infty + |f|_\infty )^{1\over 1+\alpha} \right)W_z(x) + |\psi_t|_\infty |t-\tau|.$ is a supersolution   of $1_{|f|_\infty}$ since   defining  $\lambda_2 =  {c_\psi\over \underline{ c}}  +( |\psi_t|_\infty + |f|_\infty )^{1\over 1+\alpha} $,   one has $\lambda_2 > \lambda_1 = ( |\psi_t|_\infty + |f|_\infty )^{1\over 1+\alpha}$ and then 
  \begin{eqnarray*}
   -F(x, \lambda_2 \nabla W_z,&& \lambda_2 D^2 W_z)-h(x)\cdot \lambda_2\nabla W_z |\lambda_2 \nabla W_z|^\alpha\\
   & =& -\left({ \lambda_2 \over \lambda_1}\right)^{1+\alpha } \left( F(x, \lambda_1 \nabla W_z, \lambda_1 D^2 W_z)+h(x,t)\cdot \nabla (\lambda_1 W_z )|\nabla(\lambda_1  W_z)|^\alpha\right)\\
   &\geq &-F(x, \lambda_1 DW_z, \lambda_1 D^2W_z )-h(x,t)\cdot \nabla(\lambda_1  W_z) |\nabla(\lambda_1  W_z)|^\alpha\\
   &\geq &|f|_\infty + |\psi_t|_\infty 
   \end{eqnarray*}
   Moreover in the viscosity sense , $\partial_t (|t-\tau|)\geq -1$.  This implies that 
   all the functions in the infimum are supersolutions of $1_{\{|f|_\infty\}}$.  Acting as in the proof of proposition 3  in section  4, one can prove that $W_1$ being the infimum of  supersolutions is a supersolution. 
  
We prove that $W_1$ satisfies the boundary condition on the lateral boundary  
  $W_1(x,t)  : =  \psi(x, t)$  for $x\in \partial \Omega$ and $t\in ]0,T[$ .   Indeed  first taking $(x, t)$ in the infimum  one gets $W_1(x, t)\leq \psi(x, t)$. On  the other hand  
   for all $(z, \tau) \in \partial \Omega \times ]0,T[$ 
    $ \psi(z, \tau)+ {c_\psi\over \underline{c}} \underline{c} |x-z|^\gamma + |\psi_t|_\infty |t-\tau| \geq \psi (x,t)$ 
    which implies  by considering the infimum,   the  reverse inequality. 
    
    The same arguments permit  to check that $W_1(x,0)\geq \psi (x,0)$ for all $x\in \Omega $.

   We now define  $q_1 =\sup\{2,{\alpha+2\over \alpha+1}\}$,  $q = {q_1\over \gamma}$,  $c_q = (q-1)^{q-1} + (q-1)^{1-q\over q}$.

  and also 
\begin{equation}\label{eqk2}K_2=     ({\rm diam}\ \Omega |h|_\infty + A(N+q_1-2)) ({\rm diam} \ \Omega)^{\sup (\alpha,0)},
\end{equation}
    Then, it is not difficult to see that  for any  positive constant $K_1$  and for all $y$ 
    $$(x,t)\mapsto K_1 |x-y|^{q_1}+ K_1^{1+\alpha}  K_2 t$$
    is a supersolution of $1_{\{0\}}$ and then in particular taking 
    $K_1 =  {c_\psi^q\over c_q^q\kappa^{q-1}} $
    with $c_q$ defined above,  
    for  all $\kappa\in \R^+$ and $y\in \Omega $ 
    
   $$(x,t)\mapsto {c_\psi^q\over c_q^q\kappa^{q-1}} |x-y|^{q_1} + (|f|_\infty+ |\psi_t|_\infty )  t+
    \left({c_\psi^q\over c_q^q \kappa^{q-1}}\right)^{1+\alpha} K_2 t$$
   is a supersolution of $1_{\{|f|_\infty\}}$.

Then if we define

   $$W_2(x,t) : = \inf_{y\in \Omega, \kappa\in \R^+} \{ \psi(y,0) + \kappa + {c_\psi^q\over c_q^q \kappa^{q-1}} |x-y|^{q_1} + (|f|_\infty+|\psi_t|_\infty)  t+ \left({c_\psi^q\over  c_q^q \kappa^{q-1}}\right)^{1+\alpha} K_2 t\}, $$
   
  $W_2$ being the infimum of supersolutions of $1_{|f|_\infty}$,  it is a supersolution of $1_{|f|_\infty}$. 
   
   We need to check that $W_2(x,0) = \psi (x)$.
  On  one hand, by taking $y= x$ in the infimum  and $t = 0$ one gets 
$W_2(x,t) \leq \kappa + \psi (x,0)$
 for all $\kappa$ and     on the second hand,    we  use  the  identity  for $q>1$, and for any positive number $P$ 
   \begin{equation}\label{eqkappa}\inf_{\kappa\in \R^+} \{\kappa+ {  P\over  c_q^q\kappa^{q-1}} \} = P^{1\over q} 
   \end{equation}
   
   that  we  apply  here  with $P = c_\psi^q |x-y|^{q_1}$. It gives 
   \begin{eqnarray*}
   W_2(x,0)&= & \inf_{y\in \Omega, \kappa\in \R^+} \{ \psi(y,0) + \kappa + {c_\psi^q\over c_q^q \kappa^{q-1}} |x-y|^{q_1} \}\\
   &=&  \inf_{y\in \Omega} \{ \psi (y,0) + c_\psi |x-y|^\gamma\}\\
   &\geq& \psi (x,0). 
   \end{eqnarray*}
   
   We need also to check that $W_2(x,t) \geq \psi (x,t)$ when $x\in \partial \Omega$. 
   
   For that aim we use  for all $x\in \Omega $
   $$ W_2(x,t) \geq  \inf_{y\in \Omega, \kappa\in \R^+} \{ \psi(y,0) + \kappa + {c_\psi^q\over c_q^q \kappa^{q-1}} |x-y|^{q_1} \}+ |\psi_t|_\infty |t| \geq \psi (x,0) + |\psi_t|_\infty |t| \geq \psi (x,t).$$
Moreover since $W_2$ is an infimum of continuous function it is upper semicontinuous and then  for all $x\in \partial \Omega$ and for all $t \in ]0,T[$
$$ W_2(x, t) \geq  \limsup_{x_n\in \Omega, x_n\rightarrow x}W_2(x_n, t)\geq  \lim_{x_n\in \Omega, x_n\rightarrow x}\psi (x_n, t) = \psi (x,t)$$

              We now define 
        $$W(x,t) = \inf (W_1 (x,t), W_2(x,t))$$
      Then $W$ is a supersolution of $1_{\{|f|_\infty, \psi\}}$

Similarly one can define a sub-solution  : 
$$V(x, t) = \sup (V_1(x, t), V_2(x, t))$$
with 
  $$V_1(x,t) := \sup_{(z, \tau)\in \partial \Omega \times ]0,T[,   \kappa\in \R^+} \{ \psi (z, \tau)-  \left({c_\psi\over \underline{c}} +( |\psi_t|_\infty + |f|_\infty )^{1\over 1+\alpha}\right) W_z(x) - |\psi_t|_\infty |t-\tau|\}.$$
  
and 

   $$V_2(x,t) = \sup_{y\in \Omega, \kappa\in \R^+} \{ \psi(y,0) 
  - \kappa - {c_\psi^q\over c_q^q \kappa^{q-1}} |x-y|^{q_1} - (|f|_\infty+|\psi_t|_\infty)   t- \left({c_\psi^q\over c_q^q \kappa^{q-1}}\right)^{1+\alpha} K_2 t\}$$
and $K_2$ has been defined before.Then $V$ is a sub-solution of $1_{\{-|f|_\infty, \psi\}}$. This ends the proof of proposition \ref{sub}.

    \bigskip
     
   Moreover by the  comparison principle  in theorem  \ref{thcomp}
   $$V \leq W.$$ 

\section{Existence and regularity.}

In this section, we  first prove, via Perron's method  and with the aid of  the sub and supersolutions just defined,  that there exists $u$ a unique  continuous solution of
$$\left\{\begin{array}{lc}
u_t-F (x, \nabla u, D^2 u)-h(x,t)\cdot \nabla u|\nabla u|^\alpha  = f & \mbox{in}\ Q_T\\
u=\psi(x,t)& \mbox{on}\ \partial Q_T.
\end{array}
\right.
$$
Next we prove some H\"older's estimates on this solution.

 We consider $V$ and $W$ as before,  
 $V\leq W$,  and $V$ is a subsolution, $W$ is a supersolution. 
Let
 $$E = \{ u, {\rm subsolution} \ {\rm of} \ 1_{\{ f,\psi\}}\ ,  \ V\leq u\leq W\}.$$
 Using  Perron's method adapted to our context we need to  prove that for
  $u = : \sup E$,   the lower semi-continuous enveloppe $ u_{\star}$ is a super solution of $(1)_{ f,\psi}$ , while $u ^{\star}$ is a  sub-solution. This can be done using the following proposition :
   \begin{prop}
Suppose that $\Omega$ is some open set in $\R^N$. 
 Suppose that
$u_n$ is some   locally uniformly bounded  sequence of   sub-solutions  
 for 
 $$(u_n)_t-F(x, \nabla u_n , D^2 u_n)-h(x)\cdot \nabla u_n |\nabla u_n|^\alpha   \leq f.$$
Let $\bar u$ be defined as 
 $$\bar u(\bar x, \bar t) = \limsup_{r\rightarrow 0}\{ u_n(y,s), n\geq {1\over r}, |t-s|+  |y-x|\leq r\}$$

 Suppose that $f$ is upper semicontinuous. Then $\bar u$ is a sub-solution . 
\end{prop}

Proof

$\bar u$ is upper semicontinuous by construction. 
  
   We assume that we are in the "bad " case, ie that $(\bar x , \bar t) $ is such that 
there exists $\varphi  $ which depends only on $t$,  such that $\varphi (\bar t)=0$,  and  for some $\delta_1$, 
   $\sup_{t\in B(\bar t, \delta_1)} (\bar u  (x,t)-\varphi(t)) = \bar u(\bar x, \bar t) $, 
    with  for some $\delta$, 
    $x\mapsto  \sup_{t\in B(\bar t, \delta_1)} (\bar u (x,t)-\varphi(t)) $ is  constant on $B(\bar x, \delta)$.   Then 
    $\max_{x\in B(0, \delta), t\in B(\bar t, \delta_1)} (\bar u (x,t)-\varphi(t)) =\bar  u(\bar x, \bar t)$.

    Let $k  > \sup (2, {\alpha+2\over \alpha+1})$. 
    
     We also have 
      $\sup_{x\in B(0, \delta), |t-\bar t|< \delta_1} \{\bar u(x,t)-\varphi(t))- |x-\bar x |^k- |t-\bar t|^2\}
     = \bar u(\bar x, \bar t)$ and the supremum is strict in $x$ and $t$ separately. 
     
     We now consider  
     $$\sup_{x\in B(\bar x, \delta), |t-\bar t|< \delta_1}  \{u_n^\star (x,t)-\varphi(t)-|x-\bar x|^k- |t-\bar t|^2\}
     $$ This supremum  is achieved 
      on some   $(x_n, t_n)$. 
      We begin to observe that $u_n^\star (x_n, t_n) \rightarrow \bar u(\bar x, \bar t)$. 
      Indeed by definition  of $\bar u$,  there exists $(y_n, s_n) $ which goes to $(\bar x, \bar t)$ and $u_n^\star (y_n, s_n)\rightarrow  \bar u(\bar x, \bar t)$. 
      Then 
      $u_n^\star  (x_n, t_n) -\varphi(t_n) -|x_n-\bar x|^k-|t_n-\bar t|^2 \geq u_n^\star (y_n, s_n )-\varphi(t_n)-|y_n-\bar x|^k-|s_n-\bar t|^2 \rightarrow  \bar u (\bar x, \bar t)$, which implies that 
      $\liminf u_n^\star(x_n, t_n) \geq \bar u(\bar x , \bar t)$. On the other hand,  using the definition of $\bar u$
      $$\limsup_nu_n^\star (x_n, t_n)\leq   \bar u(\bar x, \bar t).$$
    Moreover since the supremum is strict, $(x_n, t_n)\rightarrow (\bar x, \bar t)$.   
     
      If $\bar x\neq x_n$ for  an infinity of $n$, using the fact that 
      $(x,t)\mapsto \varphi(t)+ |x-\bar x|^k+ |t-\bar t|^2$ is  a test function for $u_n^\star $ on $(x_n, t_n)$ with a non zero gradient with respect to $x$ on $(x_n, t_n)$,  one gets that  for some constant $C$ 
      \begin{eqnarray*}
       \varphi^\prime (t_n) + 2(t_n-\bar t) &-& Ck^{2+\alpha} |x_n-\bar x|^{k(\alpha+1)-\alpha-2} -k^{1+\alpha} |h|_\infty 
        |x_n-\bar x|^{(k-1)(\alpha+1)}\\ 
     &\leq & \varphi^\prime (t_n) + 2(t_n-\bar t) -F( k|x_n-\bar x|^{k-2} (x_n-\bar x), D^2 (|x-\bar x|^k) (x_n)) \\
      &-&h(x_n, t_n) \cdot (x_n-\bar x) k^{1+\alpha} |x_n-\bar x|^{(k-1)(\alpha+1)-1}\\
      &\leq &f(x_n, t_n)
      \end{eqnarray*}
      This gives the result by passing to the limit since $k> {\alpha+2\over \alpha+1}$ and $f$ is upper semicontinuous.
       We now suppose that $x_n=  \bar x $ for all $n$ large enough. Then using  lemma 1  in its form for sub- solutions one gets that 
       $$\varphi^\prime (t_n) + 2(t_n-\bar t) -0 \leq f(\bar x, t_n).$$
        Once more by passing to the limit and using the upper semi continuity of $f$  we get the desired result.   
        
      When we are not in the "bad case", one can argue as in \cite{I} and \cite{BD2},  Proposition 5.2,  so we finally get  that $\bar u$ is a supersolution.

  By the comparison principle Theorem \ref{thcomp},  we get that ${u}_{\star}\geq u^{\star}$
  hence the function $ u$ is continuous and it is the required solution.
 We also know that it is unique, again by the comparison principle.

   \bigskip
   
We now prove  some  H\"older's estimate : 
    
    \begin{theo}\label{thhold}

     Let $u$ be  the  solution of $1_{ \{f,\psi\}}$. Suppose that $f$ is continuous,  bounded on $Q_T$,   and H\"older's continuous  of exponent  $\gamma_f$ with respect to $t$, that $\psi$ is  H\"older's continuous with exponent $\gamma$  with respect to $x$ and Lipschitzian in $t$. 
      Then  there exists  some constant $c$,  such that 
      for all $(x,t), (y,s)$ in $Q_T^2$,  and for $q= {q_1\over \gamma} =\sup \left( {\alpha+2\over \gamma (\alpha+1)}, {2\over \gamma}\right)$ $\gamma^\star = \inf (\gamma_f, {1\over q(\alpha+1)-\alpha })$
         $$|u(x,t)-u(y,s)|\leq c(|x-y|^{\gamma}+ |t-s|^{\gamma^\star}).$$
       \end{theo}
\begin{cor}
       Suppose that $(f_n)$  is a sequence of uniformly  bounded  functions, continuous w.r.t.   $x$ and  uniformly H\"olderian in $t$, and  $(\psi_n)$ is uniformly Holder's continuous in $x$ and uniformly Lipshitzian in $t$,  then the sequence $(u_n)$ of  solutions  of $1_{\{f_n, \psi_n\}}$ is uniformly Holder's continuous and bounded. 
 \end{cor}
      In order to prove Theorem    \ref{thhold} we give two preliminary results, which establish  some H\"older's estimates on the bottom and on the lateral boundary of $Q_T$. 
      
\begin{prop}\label{Lipt1}
Let $Q_T = \Omega \times ]0,T[$. 

Let $\psi$ be an H\"older function  with exponent $\gamma$ in $x$ and Lipshitzian in $t$ on $\partial  Q_T$, let $f$ be continuous on $\overline{Q_T}$   and let $u$ be the solution of 
$$\left\{\begin{array}{lc}
\partial_t{u}=F_{ }(x, \nabla u, D^2 u )+ h(x,t)\cdot\nabla u|\nabla u|^\alpha +  f(x, t)&{\rm in} \ Q_T\\
u(x,t)= \psi (x,t)&\ \ {\rm on} \ \ (\partial \Omega\times ]0,T[) \ \cup (\Omega \times \{0\})
\end{array}\right.
$$
Then  there exists  some constant $C_2$  such that,  for all $(x,t)\in \Omega\times ]0, T[ $ ,
$$|u(x,t)-\psi (x,0)|\leq  C_2t^{1\over q(\alpha+1)-\alpha }$$
(We recall that $q = {\sup (2, {\alpha+2\over \alpha+1})\over \gamma}$). 
\end{prop}
{\em Proof.}

By the comparison principle in theorem \ref{thcomp} one has 
\begin{eqnarray*}
u (x,t)&\leq &W(x,t)\\
&\leq & W_2(x,t) \\
&\leq &\psi (x,0) + \inf _{\kappa\in \R^+}\left(\kappa  +\left( {c_\psi^q\over c_q^q  \kappa^{q-1}}  \right)^{1+\alpha} K_2 t \right)+ (|f|_\infty + |\psi_t|_\infty )t\\
&=& \psi (x,0) + C t^{1\over (q-1)(1+\alpha)+1} +( |f|_\infty+|\psi_t|_\infty ) t
\end{eqnarray*}
 for some constant $C$  which depends on  $(c_\psi, A, a, q_1, \gamma)$, computed with the aid of (\ref{eqkappa}) replacing $q$ by $(q-1)(\alpha+1)+1$. 
  
This  yields the result.  
The  symmetric lower bound  is obtained by considering $V$ instead of $W$ and proceeding similarly.

As a consequence one has  the following

 \begin{prop}\label{Lipt}
 
 We assume here that $f$ is  continuous on $\overline{Q_T}$, H\"older    with respect to $t$,  with some exponent $\gamma_f$.   Let $u$ be a solution of $1_{\{ f,\psi\}}$. Then there exists $C_2$   depending   on the H\"older's constant $c_\psi$  and $c_f$  of $\psi$  and   $f$ respectively ,  such that  for all $x\in \Omega$  and for all $(t, s) \in ]0, T[^2$,
$$|u(x,t+s)-u(x,t)|\leq C_2s^{\gamma^\star}.$$
where  $\gamma^\star = \inf ({1\over q(\alpha+1)-\alpha},  \gamma_f)$, $q = {\sup (2, {\alpha+2\over \alpha+1})\over \gamma}= {q_1\over \gamma}$.  

\end{prop}
{\em Proof of Proposition \ref{Lipt}}: 
Let  $c_f$ be such that
$$|f(x,t+s)-f(x,t) |\leq c_f s^{\gamma_f}.$$
We define for $s
$ fixed in ]0, T[ 
\begin{eqnarray*}
v (x,t) &=& u (x, t+s) + t\ c_f s^{\gamma_f}+ \sup _{(x,t)\in \partial \Omega\times ]0, T-s[} |\psi (x, t+s)+ c_f ts^{\gamma_f}-\psi (x,t)|\\
&+&\sup_{x\in \Omega}|u (x,s)-\psi (x,0)|
\end{eqnarray*}

Then 
$v$ satisfies  on $\Omega \times ]0, T-s[$
$$\partial_t v -F (x, \nabla v, D^2 v) -h(x,t)\cdot \nabla v|\nabla v|^\alpha = f(x, t+s)+ c_f s^{\gamma_f} \geq f(x,t)$$Since $u$ satisfies the opposite inequality on the same open set,  and by construction $v(x,t)\geq u(x,t)$ on $\partial Q_T$, one has by theorem  \ref{thcomp}  

$$u(x,t) -v (x,t)\leq 0, $$
which gives the result,  redefining  $C_2= 2T ^{1+ \gamma_f-\gamma^\star}+ |\psi_t|_\infty T^{1-\gamma^\star} + C_2 T^{{1\over q(\alpha+1)-\alpha }-\gamma^\star}$  For the reverse inequality, one uses  fro $s$ fixed 
\begin{eqnarray*}
v(x,t)  &=& u (x, t+s) - t\ c_f s^{\gamma_f}- \sup _{(x,t)\in \partial \Omega\times ]0, T-s[} |\psi (x, t+s)+ c_f ts^{\gamma_f}-\psi (x,t)|\\
&
-&\sup_{x\in \Omega}|u (x,s)-\psi (x,0)|
\end{eqnarray*}

$v$ is a sub-solution of 
$$v _t -F(\nabla v, D^2 v)-h(x,t)\cdot \nabla v|\nabla v|^\alpha  \leq f(x, t+s)-c_f s^{\gamma_f}\leq f(x, t)$$
and  
$u (x, t)$ satisfies the opposite inequality on $]0, T-s[$.  Moreover $v(x,t)\leq u(x,t)$ on $\partial Q_T$. 
Then Theorem \ref{thcomp} implies that 
$$
u (x, t+s)
\leq u (x,t)+  C_2 s^{\gamma^\star}
$$
with  $C_2$ as above. 

We now give an estimate on the lateral boundary : 

\begin{prop}\label{holx} We assume that $\psi$ is H\"older continuous  of exponent $\gamma$  with respect to $x$ and Lipschitzian with respect to $t$. Let $u $ be a solution of $1_{\{f , \psi\}}$.  Then  there exists $C_1$  such that for  all $(x,x_o)\in \Omega\times \partial \Omega$ and $t\in [0,T)$,

$$|u (x,t)-u (x_o,t)|\leq C_{1}|x-x_o|^{\gamma}.$$
\end{prop}
{\em Proof}

We use once more  the supersolution.   
Taking in the infimum defining $W$ the point $(x_o, t)$ which is on the lateral boundary, and using the properties of the  barrier,  one has 
\begin{eqnarray*}
u(x,t)&\leq& W(x,t)\\
&\leq& W_1(x,t)\\
&\leq & \psi (x_o,t)+{c_\psi\over \underline{c}} W_{x_o}(x)  +( |f|_\infty + |\psi_t|_\infty )^{1\over 1+\alpha} W_{x_o} (x)\\
&\leq & \psi (x_o, t)+({c_\psi\over \underline{c}} +  (|f|_\infty+ |\psi_t|_\infty )^{1\over 1+\alpha})\bar c |x-x_o|^\gamma.
\end{eqnarray*}
This gives the result with 
$$C_1 =\bar c ({c_\psi\over \underline{c}} +  (|f|_\infty+ |\psi_t|_\infty )^{1\over 1+\alpha})$$  One gets the lower bound by considering  $V$ instead of $W$. 
\bigskip

We now prove Theorem  \ref{thhold}. 
First observe that $u$ is bounded as soon as $f$ and $\psi$ are bounded, due to  theorem \ref{thcomp}, the inequalities 
$V\leq u\leq W, $
 and the definition of $V$ and $W$. 

In the following $\delta$ will be $< \inf (1, {1\over T})$,  and $L>1$. 

We construct a function $\Phi$ as follows:
Let $\delta$ be small enough in order that,  for $\tilde\omega$ the modulus of continuity  given in  the assumption (H3),  and $C$ being the universal constant defined in (\ref{elip1}) later,  one has 
$\tilde\omega (\delta)< {a\over 4 C}$, and $\delta |h|_\infty < {a\over C}$. We define  
  $$L= \sup  \left(C_1,\displaystyle{ \left({|f|_\infty \delta^{\alpha+1-(\alpha+2)\gamma} \over a \left(\gamma\right)^{1+\alpha} (1-\gamma)}\right)^{1\over 1+\alpha}}, {2\sup u\over \delta^{\gamma}}\right)$$
  $$M = \sup (TC_2, {2\sup u\over \delta^{\gamma^\star }})$$
   where $C_1$ is given in Proposition \ref{holx}, 
  and $C_2$ is given in Proposition \ref{Lipt}. We also define  
$$\Delta_\delta=\{((x,t),(y,s))\in Q_T^2,\ |x-y|<\delta,|t-s|< \delta  \}.$$

\noindent {\bf Claim} {\em For any} $(x,t),(y,s)\in\Delta_{\delta}$
\begin{equation}\label{ho}
\Phi(x,t,y, s)= u (x,t)-u(y,s)- L |x-y|^{\gamma }-M|t-s|^{\gamma^\star }\leq 0.
\end{equation}

Suppose for a while that the supremum of $\phi$ is positive. Then, for $\kappa$ small enough  the supremum of $\phi-{\kappa\over T-t}-{\kappa\over T-s}$ is also strictly positive. In the following we replace $\phi$ by $\phi-{\kappa\over T-t}-{\kappa\over T-s}$ .

From the choice of the constants and Propositions  \ref{Lipt}  and \ref{holx} we know that the inequality 
(\ref{ho}) with the "new " $\phi$ holds on $\partial\Delta_{\delta}$ : 

Indeed if $x\in \partial \Omega $, $y\in \Omega$,  and $(t,s)\in ]0, T[^2, |s-t|< \delta$, using 
Proposition   \ref{holx}, one has 
\begin{eqnarray*}
u(x,t)-u(y,s)&\leq & \psi(x,t)-\psi (x, s)+ u(x,s)-u(y,s)\\
&\leq & |\psi_t|_\infty |t-s| + C_1 |x-y|^\gamma
\end{eqnarray*}
 which gives the result since $M\geq C_2 \geq |\psi_t|_\infty$ and $L \geq C_1$. 
The same is true  by exchanging $x$ and $y$. 

 If $|x-y|=\delta$ or $|t-s|= \delta, $ the result holds by the choice of $L$ and $M$. 
 For $t=0$ or $s=0 $, one uses  proposition \ref{Lipt} and proposition \ref{holx} to get 
 $|u(x,t)-u(y,0)|\leq |u(x,t)-u(x,0)|+ |u(x,0)-u(y,0)|\leq c_\psi |x-y|^\gamma+ C_2 t^{\gamma^\star}$, from  which we conclude since $L > c_\psi$ and $M > C_2 T$. 
  
Finally    the  supremum cannot be achieved  on $t = T$  or $s=T$ since in that case the function is $-\infty$.

Suppose by contradiction that
$$\sup_{(x,t), (y,s)\in Q_T^2}\Phi(x,t,y,s)>0.$$

Then 
 for $n>0$ large  enough 
 $$\Phi_n  (x,t,y,s) =u (x,t)-u(y,s)- L |x-y|^{\gamma}-M(|t-s|^2+ n^{-2})^{\gamma^\star\over 2}-{\kappa\over T-t}-{\kappa\over T-s}$$
 has also a supremum $>0$, and it cannot be achieved on the boundary, by the previous considerations. We denote for simplicity  by $(\bar x_n, \bar t_n), (\bar y_n, \bar s_n)$ a couple   inside $\Delta_\delta$ on which the supremum  of $\psi_n$ is achieved. In the following we fix $n$ large enough and drop the indexes $n$ for simplicity. 
  
Suppose that $\bar x = \bar y$. Then one would have 
$$u(\bar x,t)-u(\bar x,s)\geq  {M((t-s)^2+{1\over n^2})^{\gamma^\star\over 2}  }, $$
which contradicts  proposition \ref{Lipt} and the choice of $M$. Hence  $\bar x\ne \bar y$ and 
  using Ishiis' lemma (see also lemma 2.1 in \cite{BD1}), 
there exists 
$X\in S$ and $Y$ in $S$ such that:

$$\left(M  \gamma^\star {\bar t-\bar s\over ((\bar t-\bar s)^2+{1\over n^2} )^{1-{\gamma^\star\over 2}}}+ {\kappa\over (T-\bar t)^2} ,  \gamma L(\bar x-\bar y)|\bar x-\bar y|^{\gamma-2}
, X\right)\in J^{2,+ }u(\bar x,\bar t)$$
$$\left( M  \gamma^\star {\bar t-\bar s\over ((\bar t-\bar s)^2+{1\over n^2} )^{1-{\gamma^\star\over 2}}}-{\kappa\over (T-\bar s)^2},\gamma L(\bar x-\bar y)|\bar x-\bar y|^{\gamma-2}
, -Y\right)\in J^{2,-}u(\bar y,\bar s)$$
with 
$$\left(\begin{array}{cc}
X&0\\
0&Y
\end{array}\right) \leq \left(\begin{array}{cc}
B &-B\\
-B&B
\end{array}\right) $$
and  $B  = L\gamma|x-y|^{\gamma-2} (I+ (\gamma-2){(x-y)\otimes (x-y)\over |x-y|^2}) = D^2 (|X|^{\gamma})(x-y)$.

 We need a 
more precise estimate, as in
\cite{IL}. For that aim let  $P$ be defined as :

$$0\leq P : = {(\bar x-\bar y\otimes \bar x-\bar y)\over |\bar x-\bar y|^2}\leq I.$$

Using $-(X+Y)\geq 0$,  $(I-P)\geq 0$ and the properties of the
symmetric matrices one has 
$$tr(X+Y)\leq tr(P(X+Y)).$$ Remarking in addition that
$X+Y\leq 4B$, one  sees that
$tr(X+Y)\leq tr(P(X+Y))\leq 4tr(PB)$. But $tr (PB)=\gamma L(\gamma-1)|
\bar x-\bar y|^{\gamma-2}<0$, hence

\begin{equation}\label{elip}
|tr (X+Y)|\geq 4\gamma L(1-\gamma)| \bar x-\bar y|^{\gamma-2}.
\end{equation}

Furthermore  by Lemma III.1 of \cite{IL} there exists a universal constant $C$ such that 

\begin{equation}\label{elip1}
|X|, |Y|\leq C (|tr(X+Y)|+ |B|^{1\over 2} |tr(X+Y)|^{1\over
2})\leq C|tr(X+Y)|
\end{equation}
since $|B|$ and $|tr(X+Y)|$ are of the same order.  This constant is the constant used for the choice of $L$ chosen  at the beginning of the proof.

Using the  fact that $u$ is both a sub- and a supersolution  we get

\begin{eqnarray*}
f(\bar x, \bar t) &\geq& {M\gamma^\star}\left({\bar t-\bar s\over ((\bar t-\bar s)^2+{1\over n^2})^{1-{\gamma^\star\over 2}}}\right)+ {\kappa\over (T-\bar t)^2}\\
&
- & F (\bar x, \gamma L(\bar x-\bar y)|\bar x-\bar y|^{\gamma-2}, X)\\
&-&L^{1+\alpha}\gamma^{1+\alpha}  h(\bar x, \bar t) \cdot (\bar x-\bar y)|\bar x-\bar y|^{(\gamma -1)(\alpha+1)-1}\\
 &\geq&   {M\gamma^\star}\left({\bar t-\bar s\over ((\bar t-\bar s)^2+{1\over n^2})^{1-{\gamma^\star\over 2}}}\right)- {\kappa\over (T-\bar s)^2}- F (\bar y, (\gamma L(\bar x-\bar y)|\bar x-\bar y|^{\gamma-2}, -Y)\\
 &-&L^{1+\alpha}\gamma^{1+\alpha} h(\bar y, \bar s) \cdot (\bar x-\bar y)|\bar x-\bar y|^{(\gamma -1)(\alpha+1)-1}-\tilde\omega(|\bar x-\bar y|)
(\gamma L  |\bar x-\bar y|^{\gamma-1})^\alpha|X|  \\
&-&L^{1+\alpha}|h|_\infty\gamma^{1+\alpha} |\bar x-\bar y|^{(\gamma)(\alpha+1)}
+(\gamma L |\bar x-\bar y|^{\gamma-1})^\alpha a|tr(X+Y)|\\
&\geq & f(\bar y, \bar s) +4\gamma^{1+\alpha} L^{1+\alpha } (1-\gamma) |\bar x-\bar y|^{\gamma-2+ (\gamma-1)(\alpha+1)}(a-{\tilde\omega\over C} (|\bar x-\bar y|) -{|h|_\infty\over 4C}| \bar x-\bar y|)
\end{eqnarray*}
which is a contradiction with the assumptions on $L$. 
We have obtained that 
$$u(x,t)-u(y,s) \leq L |x-y|^{\gamma} +M {|t-s|^{\gamma^\star}\over T-t}.$$
This ends the proof.

 \section{Maximal solutions on $\Omega \times \R^+$}
 
 In this section we prove the existence of solutions  on $\Omega\times \R^+$.   For this we prove some property of solutions when $t\rightarrow T$ and we use Zorn's lemma.  
  
  \begin{prop}
  We suppose that $f$ is continuous and bounded on $ \Omega \times \R^+$. 
   Suppose that $u$ is a supersolution of $1_{\{f, \psi\}}$on $Q_T$, lower semicontinuous,  and  we define 
   $$u(x,T) = \liminf_{|z-x|+ |t-T|\leq r} u(z,t). $$
   
   Then   $u$ being extended in that kind  is a supersolution on  $\Omega \times ]0,T]$. 
   
   In the same manner  if $v$ is a upper semicontinuous  sub-solution,   we define 
   $$v(x,T) = \limsup_{|z-x|+ |t-T|\leq r} v(z,t).$$
   Then $v$ being extended in that kind  is  a sub-solution on $\Omega \times ]0,T]$.

   \end{prop}
   Proof 
   
   We follow partly the process employed in \cite{OS}. 
   
      Let $u$ be a supersolution and let $\varphi$ be a ${\cal C}^2$ function such that 
   $$(u-\varphi )(x,t) \geq (u-\varphi) (\bar x, T)$$
   for $(x,t) $ on some neighborhood $V$ of $(\bar x,T)$,  $\nabla_x \varphi (\bar x, \bar t)\neq 0$. One can assume replacing if necessary $\varphi(x,t)$ by $\varphi(x,t) - |x-\bar x|^k- |t-T|^2$ for $k>\sup (2, {\alpha+2\over \alpha+1})$,   that the infimum  of $(u-\varphi)$ is strict on $(\bar x, T)$.

    Then for $n$ large enough 
 $$\inf_{(x,t)\in V} \left(u(x,t)-\varphi(x,t)+ {1\over n(T-t)}\right)$$
 is achieved on $(y_n, t_n)$ with $(y_n,t_n)\rightarrow (\bar x, T)$. 
     
       Indeed we prove first that $$\lim_{n\rightarrow +\infty} \inf_{(x,t)\in V} \left(u(x,t)-\varphi(x,t)+ {1\over n(T-t)}\right) = \inf_{(x,t)\in V} (u-\varphi)(x,t).$$
     We already have 
     $$\inf_{(x,t)\in V} \left(u(x,t)-\varphi(x,t)+ {1\over n(T-t)}\right) \geq  \inf (u-\varphi)(x,t).$$
     
     For the reverse inequality let $\epsilon $ be given and $(x_\epsilon, t_\epsilon )$ in $Q_T$ with 
$$(u-\varphi) (x_\epsilon, t_\epsilon)\leq \inf_{(x,t)\in V}(u-\varphi) (x,t)+\epsilon$$ then  for
 $n (T-t_\epsilon) > {1\over \epsilon}$
       $$  (u-\varphi) (x_\epsilon, t_\epsilon)+ {1\over n(T-t_\epsilon)} \leq ( u-\varphi)(x_\epsilon, t_\epsilon)+ 2\epsilon\leq \inf_{(x,t)\in V} (u-\varphi) + 2\epsilon.$$
       $\epsilon $ being arbitrary, one gets the result.

      Now the function $u-\varphi+ {1\over n(T-t)}$ being lower semi-continuous  the infimum is achieved on some $(y_n, t_n)$. 
        By the previous considerations 
        $$\inf_{(x,t)\in V} (u-\varphi)(x,t) \leq (u-\varphi)(y_n, t_n) + {1\over n(T-t_n)} \rightarrow (u-\varphi) (\bar x, T)$$
 This implies in particular that 
        $$(u-\varphi)(y_n, t_n) \rightarrow (u-\varphi) (\bar x, T)$$ and since the infimum  of $u-\varphi$ is strict,  
 $(y_n, t_n)\rightarrow (\bar x, T)$.  Let us note that $t_n$ does not go to $T$ too quickly, since 
 $n(T-t_n)\rightarrow +\infty$.

  Let $\varphi_n = \varphi (x,t)-{1\over n(T-t)}$,  since $\varphi$ is ${\cal C}^1$,  for $n$ large enough,  $\nabla_x \varphi_n (y_n, t_n)\neq 0$, and since $\varphi_n$ achieves $u$ by below on $(y_n , t_n)$, 
$${d\over dt} \varphi_n(y_n, t_n)-F(y_n, \nabla \varphi(y_n, t_n), D^2\varphi(y_n, t_n))- h(y_n)\cdot \nabla \varphi(y_n)|\nabla \varphi (y_n)|^\alpha   \geq f(y_n, t_n), $$
hence 
\begin{eqnarray*}
{d\over dt} \varphi (y_n, t_n)&-&F(y_n, \nabla \varphi (y_n, t_n), D^2\varphi(y_n, t_n))- h(y_n)\cdot \nabla \varphi(y_n)|\nabla \varphi (y_n)|^\alpha\\
&  \geq &f(y_n, t_n) + {1\over n(T-t)^2}\\
&\geq& f(y_n, t_n), 
\end{eqnarray*}

and passing to the limit one gets that 
$${d\over dt} \varphi (\bar x, T)-F(\bar x, \nabla \varphi , D^2\varphi) (\bar x, T) - h(\bar x )\cdot \nabla \varphi(\bar x,T)|\nabla \varphi (\bar x,T)|^\alpha \geq f (\bar x, T).$$
This ends the case $\nabla_x \varphi(\bar x, T) \neq 0$.

We now assume that there exists some ${\cal C}^1$ function $\varphi$ which depends only on $t$,  and some $\delta_1>0$ such    that  $u(\bar x,T)-\varphi(T)  = \inf_{|t-\bar t|< \delta_1} ( u(x,t)-\varphi (t))$ and   $\inf_{|t-\bar t|< \delta_1} \{ u(x,t)-\varphi (t)\}$ is constant in a neighborhood $B(\bar x, \delta)$ of $\bar x$. Then one also has 
$$\inf  _{x\in B(\bar x, \delta), |t-\bar t|< \delta_1}\{ u(x,t)-\varphi (t)+|x-\bar x|^k+|t- T|^2\}= u(\bar x, T)-\varphi(T)$$
Defining  $\varphi_n(t) = \varphi (t)-|x-\bar x|^k-|t-T|^2-{1\over n(T-t)}$   one gets also that there exists   $(x_n, t_n)$ which converges to $(\bar x, T)$ and 
 $(x_n, t_n)$ is a local  minimum for $u-\varphi_n$.

-  Either $x_n = \bar x$ for all  $n$ large enough, then  using  lemma  1 one  gets 
  $$\partial_t \varphi (t_n) -2(t_n- T) - {1\over n(T-t_n)^2}\geq f(\bar x, t_n).$$
which yields  the result by passing to the limit.

-Or  for an infinity of $n$,  $ x_n\neq \bar x$, then 
\begin{eqnarray*}
\partial_t \varphi (t_n) -2(t_n- T) - {1\over n(T-t_n)^2} &-&F(-k|x_n-\bar x|^{k-2} (x_n-\bar x), -D^2(|\bar x-x|^k)(x_n)) \\
&+& k^{1+\alpha } h(x_n)\cdot (x_n-\bar x) |x_n-\bar x|^{(k-1)(\alpha-1)-1}\\
& \geq& f(x_n, t_n).
\end{eqnarray*}

Since $|\bar x-x_n|$  and $|t_n-T|$ tend  to zero when $n$ goes to infinity,   and $k > {\alpha+2\over \alpha+1}$, 
 one gets by passing to the limit  that 
$$\varphi^\prime (T)\geq f(\bar x, T).$$

   \bigskip
   
   We can now use Zorn's axiom to get the existence of maximal solutions  on $\Omega \times \R^+$
 for the problem $1_{\{f, \psi\}}$. 
 
 Moreover one can prove using uniform Holder's estimates that the solutions are locally Holder's on $\Omega \times \R^+$.  We do not give the proof which uses  both some arguments in  the Holder's proof for $\Omega \times ]0,T[$ and some arguments specific to the non bounded cases, as those used for the case of $\R^N$ in theorem 3 later.

 \section{The case $\R^N \times ]0,T[$}
 
 For completeness sake we are going to prove   some existence's  result  for the equation in $\R^N\times ]0,T[$ when $f$ is uniformly continuous and bounded on $\R^N\times ]0,T[$, Holder's continuous  in $t$, uniformly w.r.t. $x$,  and $\psi$ is  H\"olderian  for some exponent $\gamma_\psi$ and uniformly bounded on $\R^N$.  We assume in addition that $F$  satisfies the uniform Lipschitz condition :
 
 $(H6)$  There exists some constant $C$ such that for all $p\neq 0$,  for all $X$ and for all $q,$ such that $|q|< {|p|\over 2}$, one has 
 $$|F(x, p+q, X)-F(x, p, X)|\leq C |p|^{\alpha -1} |q| |X|$$
   We prove the existence of viscosity solutions of 
  $$\left\{\begin{array}{lc}
  u_t-F (x, \nabla u, D^2 u)-h(x,t)\cdot \nabla u|\nabla u|^\alpha  = f(x,t)&\ {\rm  in} \ \R^N \times ]0,T[\\
   u(x,0) = \psi (x)\ & {\rm on } \ \R^N \times \{0\}
   \end{array}\right.$$

We will construct  a supersolution and a sub-solution and use Perron's method to conclude. 
    
     To construct a supersolution, we use  the following proposition 
   \begin{prop}
  There exists   $G$,   some positive ${\cal C}^2$ function on $[0,\infty[$,  and some constant $B$ such that $u(x)=G(|x|)$  satisfies  on $\R^N\times ]0,T[$ 
   $$F(x, \nabla u, D^2 u)+ h(x,t)\cdot \nabla u|\nabla u|^\alpha \leq B.$$
  \end{prop}
{\em    Proof} :
 If  $\alpha \geq 0$ let $G$ be defined as 
   $$G(r) = \left\{ \begin{array}{lc}
   r^2 &{\rm if} \ r<1\\
  (r-1) (3-{1\over r})+1 & {\rm if} \ r\geq 1. 
  \end{array}\right.$$
  
   In the case    where $\alpha <0$, we recall that $q_1 = {\alpha+2\over \alpha+1}$,  $q = {q_1\over \gamma_\psi}$, 
    and define 
   $$G(r) =\left\{ \begin{array}{lc}
    r^{q_1}& \ {\rm if}  \ r<1\\
    {q_1(1+q_1) r\over 2} + {q_1(q_1-1)\over 2r} +  1-q_1^2\ & {\rm if} \ r>1. 
    \end{array}
    \right.$$
 With this choice of $G$ by a tedious but straithforward computation there exists some constant $B$   such that  for $ u(x) = G(|x|)$
     $$F(x, \nabla u, D^2 u)+ h(x,t)\cdot \nabla u|\nabla u|^\alpha \leq B.$$

\bigskip
 We now define on the model of $W_2$ in section 3,   
   $$W(x,t) = \inf_{y\in \R^N, \kappa\in \R^+} \{ \psi(y) + \kappa + {(c_\psi + 2| \psi|_\infty )^q \over c_q^q\kappa^{q-1}} G(|y-x|) + |f|_\infty  t+ \left({(c_\psi + 2| \psi|_\infty )^q \over c_q^q \kappa^{q-1}}\right)^{1+\alpha} B t\}$$
  
 Then $W$  is an  infimum of supersolutions for $1_{\{|f|_\infty\}}$

Moreover  
$$ W(x,0)= \inf_{\{|y-x|<1, \kappa\in \R^+\}} (\psi (y) +\kappa + {(c_\psi + 2| \psi]_\infty )^q \over c_q^q\kappa^{q-1}} G(|x-y|) \geq 
\psi(y) + c_\psi |x-y|^{\gamma_\psi}\geq  \psi(x)$$

and also  using $G(r)\geq r$ for $t\geq 1$ 
 $$ \inf_{  |y-x|>1} \{ \psi (y) + (c_\psi+ 2|\psi|_\infty ) |y-x|^{1\over q} \}\geq \psi (y) + 2 | \psi|_\infty  \geq \psi(x).$$
This implies that 
$W(x,0) \geq \psi (x)$. 
 Moreover taking $y = x$ in the infimum, one gets 
 $$W(x,0) \leq \kappa + \psi (x), $$
 for all $\kappa$.  We have obtained that $W(x,0) = \psi (x)$. 
 We now observe that $W$ is uniformly bounded, indeed 
\begin{eqnarray*}
 W(x,t)&\leq& \inf _{|x|<1}\{\psi (x) + \kappa+  |f|_\infty  t+\left({(c_\psi + 2| \psi|_\infty )^q \over c_q^q \kappa^{q-1}}\right)^{1+\alpha} B t\}\\
  &\leq & \psi (x) + c t^{1\over q(\alpha+1)-\alpha } \\
  &\leq & |\psi|_\infty + c T^{1\over q(\alpha+1)-\alpha } 
  \end{eqnarray*}
  We do not give  explicitely $c$ which can be 
   computed using (\ref{eqkappa}), replacing $q$ by $(q-1)(\alpha+1)+1$. 
  Moreover  there exists $c_1$ and $c_2$ such that 
    $$W (x,t)\leq \psi (y)+ c_1(|x-y|^{\gamma_\psi} )+ c_2 t^{1\over q(\alpha+1)-\alpha }$$
    Indeed 
    $$W(x,t) \leq \psi (x) + c t^{1\over q(\alpha+1)-\alpha } \leq \psi (y) + c_\psi |x-y|^{\gamma_\psi} +  c t^{1\over q(\alpha+1)-\alpha  }.  $$

Let us note that 
$$V(x,t) = \sup_{y\in \R^N, \kappa\in \R^+} \{ \psi(y,0) - \kappa - {c_\psi^q\over c_q^q \kappa^{q-1}} G(|y-x|) - |f|_\infty  t- \left({c_\psi^q\over c_q^q \kappa^{q-1}}\right)^{1+\alpha} B t\}$$

 with  $B$ as before , 
is a sub-solution of $1_{\{-|f|_\infty, \psi\}}$. Moreover $V$ is bounded and satisfies for some constants $c_1$ and $c_2$
  $$V (x,t)\geq \psi (y)- c_1|x-y|^\gamma - c_2 t^{1\over q(\alpha+1)-\alpha }. $$
                          
        \bigskip
        A first crucial step  for the existence of solutions for the Dirichlet problem is some  comparison theorem on $\R^N\times ]0,T[$.  This  will also  permit to get the uniqueness and later the regularity of the solutions. 

  \begin{theo}\label{compRN}
Suppose that $f$ and $g$  are  uniformly continuous  and bounded  and $f\geq g$.
   Suppose that $u$ and $v$ are respectively  uppersemicontinuous and lower semicontinuous sub-and supersolutions of 
   $$ u_t-F(x, \nabla u, D^2 u)-h(x,t)\cdot \nabla u|\nabla u|^\alpha  \leq g(x,t) \ {\rm in} \  \R^N \times ]0,T[\\
               $$
  $$v_t-F(x, \nabla v , D^2 v)-b(x,t)\cdot \nabla v|\nabla v|^\alpha \geq f(x,t)\ {\rm in}\ \R^N\times ]0,T[$$
with $u(x,0)\leq v (x,0)$, $x\mapsto u(x, 0)$  and $x\mapsto v(x,0)$ being Holder's continuous and bounded.  Suppose in addition that  there exist some constant $c_1$,   such that  for all $x$, $y$ in $\R^N$
 $$u (x,t)\leq u (y,0)+ c_1(|x-y|+1)  $$
 $$ v (x,t) \geq v  (y,0) -c_1(|x-y|+1) $$
 Then $u(x,t)\leq v(x,t)$. 
 \end{theo}
 
 We postpone the proof of theorem \ref{compRN} and  derive from it some consequences.  
 
 First the estimates on $V$ and $W$ imply that 
 $V \leq W$. Then   using  Perron's method  in section 4, which proof  does not use the boundedness of $\Omega$,  we obtain that there exists a solution of $1_{\{f, \psi\}}$ on $\R^N\times ]0,T[$, in the sense that $u^\star$ is a sub-solution and $u_\star$ is  a supersolution.  We now use the fact that $V \leq u^\star$ and $u_\star \leq W $ to derive that 
           there exist $c_1$ and $c_2$ such that 
    $$u_\star (x,t)\leq \psi (y)+ c_1|x-y|^{\gamma_\psi} + c_2 t^{1\over q(\alpha+1)-\alpha }$$
    $$ u^\star (x,t) \geq \psi   (y) -c_1|x-y|^{\gamma_\psi} -c_2 t^{1\over q(\alpha+1)-\alpha}$$
    From these estimates,  using theorem \ref{compRN} one gets that 
    $u_\star \geq u^\star$, hence $u$ is continuous.   Applying once more theorem \ref{compRN} one gets that the solution is unique.

 {\em Proof  of theorem \ref{compRN}}
                
One can replace $v$ by $(v)_\kappa  = v + {\kappa \over T-t}$. Then $v_\kappa $ is a strict supersolution, which is infinite on $t=T$. 
                
  We shall prove that $u\leq v_\kappa$ and next we shall let $\kappa$ go to zero. 
    In the following we drop the index $\kappa$ .

Suppose by contradiction that there exists  $(\bar x, \bar t)$ such that  $(u-v)(\bar x, \bar t)>0$. 
 Then $\bar t< T$  according to the  previous property of $v$.

                We introduce  for $j\in \N$ and for $k = \sup (3, {|\alpha|\over 3}, {\alpha+2\over \alpha+1}, {\alpha+1}, {\alpha+2\over 6}, {2(1+\alpha)\over \omega_h})$,   the function $\psi_j$ defined as 
 $$\psi_j(x,y,t,s) = u(x, t)-v(y,s) -{j|x-y|^k\over k} -{1\over  j^{3k}} |x|^k-{j\over 2} |t-s|^2$$
 Then for $j$ large enough the supremum of $\psi_j$ is still $>0$, for example as soon as
   $$j^{3k} > {|\bar x|^k\over u(\bar x, \bar t)-v(\bar x, \bar t) }$$

 In the following $C$ will denote some constant which can vary from one line to another. 
                 
We prove first that  if $\psi_j(x_j, y_j,t_j, s_j)>0$, $j|x_j-y_j|^{k}\leq C$. 
 Indeed  one has  for $j^{3k}>{|\bar x|^k\over (u(\bar x, \bar t)-v(\bar x, \bar t))}$, 
 $\psi_j(x_j, y_j, t_j, s_j) \geq 0$  and then  using   $u (x,t)\leq u (y,0)+ c_1(|x-y|+1)  $
 and $$ v (y,t) \geq v  (y,0) -c_1, $$
 one gets 
 
                 \begin{eqnarray*}
                   { j|x_j-y_j|^k \over k}&\leq& c_1(|x_j
                   -y_j|+2)\\
                   &\leq & {j|x_j-y_j|^k\over 2 k} + C,
                   \end{eqnarray*}
                    and then $j|x_j-y_j|^k$  is bounded. 
                 
In particular $|x_j-y_j|$ goes to zero.  From this one also  derives that 
 $${|x_j|^k\over j^{3k}}+ {j|t_j-s_j|^2\over 2} \leq C,$$
 and then $|x_j|\leq \sqrt{C}j^{3}$

  Moreover  using Ishii's lemma \cite{I},  (see also lemma 2.1 in \cite{BD1})  there exist $(X_j, Y_j)\in {\cal S}$ such that 
  $$ \left(j(t_j-s_j), j|x_j-y_j|^{k-2} (x_j-y_j)+ {k|x_j|^{k-2} x_j\over j^{3k}}, X_j+ {D^2(|x|^k)(x_j)\over j^{3k}}\right)\in J^{2,+} u(x_j, t_j)$$
  $$ \left(j(t_j-s_j), j|x_j-y_j|^{k-2} (x_j-y_j), -Y_j\right)\in J^{2,-} v(y_j, s_j)$$

                     Suppose that  $x_j= y_j$. We prove then that $x_j \neq 0$. If it was  the case      the function 
                     $\varphi(x,t) =  u(0, t_j)+ {j\over k} |x|^k +{|x|^{k}\over j^{3k}}    
                    +{j\over 2} (t-s_j)^2- {j\over 2} (s_j-t_j)^2$
                     would touch $u$ by above on $0$ and then using lemma \ref{lem1}
 one would obtain since $k> \sup(2, {\alpha+2\over \alpha+1} )$
 $$j(t_j-s_j)-0 \leq g(0, t_j).$$
 On the other hand since    
$ -{j\over 2} (t-t_j)^2+ {j\over 2} (s_j-t_j)^2-{j|x|^k\over k}$ touches $v$ by below on $(0, s_j)$,  using once more lemma \ref{lem1} we get 
                      $$j(t_j-s_j)-0 \geq f(0, s_j)+ {\kappa\over (T-s_j)^2}$$
                      
                     Using $|t_j-s_j|\rightarrow 0$,  the uniform continuity of $f$ and $g$ , substracting  the two inequalities and passing  to the limit we get a contradiction. 

                       We now suppose that $x_j = y_j$ and  we know that  under this assumption, $x_j \neq 0$. Then  
                       the function    $\varphi(x,t) =v(x_j, s_j)+  {j\over k} |x-x_j|^k +{|x|^{k}\over j^{3k}}    
                    +{j\over 2} (t-s_j)^2- {j\over 2} (s_j-t_j)^2$
                    achieves  $u$ by above on $x_j$,  where its  gradient is different from $0$. 
                    We then have 
                    $$j(t_j-s_j) -F({k|x_j|^{k-2} x_j\over j^{3k}}, D^2 ({|x|^k\over j^{3k}}))(x_j)-\left(h(x_j, t_j)\cdot {k|x_j|^{k-2}x_j \over j^{3k}}\right)\left({|k|x_j|^{k-1} \over j^{3k}}\right)^\alpha \leq g(x_j, t_j), $$
                    and for $v$ one uses once more  lemma \ref{lem1}, 
                    to  get that 
                     $$j(t_j-s_j)-0 \geq f(x_j, s_j)+ {\kappa\over (T-s_j)^2}, $$
                      We now use the properties of $F$ to get that 
                      \begin{eqnarray*}
                     \vert F({k|x_j|^{k-2}x_j \over j^{3k}}, D^2 {(|x|^k)\over j^{3k}}(x_j)))&+&\left(h(x_j, t_j)\cdot {k|x_j|^{k-2} x_j\over j^{3k}}\right)\left({k|x_j|^{k-1} \over j^{3k}}\right)^\alpha\vert \\
                      &\leq& C\left({ |x_j|^{k(\alpha+1)-\alpha-2}\over j^{3k(1+\alpha)}}+ { |x_j|^{(k-1)(\alpha+1)}\over j^{3k(1+\alpha)}}\right)\\
                      &\leq &C \left(
                      j^{3(k(\alpha+1)-\alpha-2)-3k(1+\alpha)}+ j^{-3(\alpha+1)}\right) \\
                      &=& o(1)
                      \end{eqnarray*}
                      
                              Finally using the fact that $|x_j-y_j|+ |t_j-s_j|$ goes to zero,   the uniform continuity of $f$ and $g$ ,  substracting  the two equations and passing  to the limit we  get a contradiction.

                                         We have obtained that $x_j \neq y_j$. 
                    
                 We now prove that 
                        $j^{2}|x_j-y_j|^{k-1} \rightarrow +\infty$. In particular this will imply that for $j$ large enough 
                        $j|x_j-y_j|^{k-2} (x_j-y_j)+ k{|x_j|^{k-2} x_j\over j^{3k}}\neq 0$.  Suppose by contradiction  that  for some constant $c>0$, 
                        $j|x_j-y_j|^{k-1} \leq c j^{-1}$
                         then 
                        $|X_j|\leq j|x_j-y_j|^{k-2} \leq (j^2|x_j-y_j|^{k-1})^{k-2\over k-1} j^{3-k\over k-1} \rightarrow 0$ and also $|X_j|+ |D^2\left({|x|^k\over j^{3k}}\right) (x_j)|\leq |X_j|+ cj^{-6}  \rightarrow 0$. Using the fact that $u$ and $v$ are respectively sub-and supersolution, one has  
                      
                         \begin{eqnarray*}
                         g(x_j, t_j)&\geq&  j(t_j-s_j) - o(1)\                  \\
                         {\rm       and}\   {\kappa\over T^2} +  f(y_j, s_j)  &\leq &j(t_j-s_j)+ o(1).
                          \end{eqnarray*}
                             Substracting the two inequalities, passing to the limit and using the properties of $f$ and $g$,  one gets a contradiction.  We have obtained that $j|x_j-y_j|^{k-1} \geq {c\over j}$ for some constant $c$.  From this one derives that $j|x_j-y_j|^{k-2} (x_j-y_j) + {k|x_j|^{k-2} x_j\over j^{3k}}\sim_{j\rightarrow +\infty}  j|x_j-y_j|^{k-2} (x_j-y_j) $. 
                                                        With the aid of this remark and  using the assumption    $(H6)$                   \begin{eqnarray*}
                         & | &F (j|x_j-y_j|^{q-2} (x_j-y_j)+ {k|x_j|^{k-2} x_j\over j^{3k}},X_j ) -F (j|x_j-y_j|^{q-2}(x_j-y_j),   X_j)|\\
                          &\leq&c j^{-3} |X_j|(j|x_j-y_j|^{k-1})^{\alpha-1}\\
                          &\leq & \left\{ \begin{array}{cc}
                          c j^{-\alpha-1}|x_j-y_j|^{k-2}   & {\rm if} \ \alpha <1\\
                            cj^{-3+ {\alpha+1\over k}} (j|x_j-y_j|^k)^{\alpha-{\alpha+1\over k}}  &{\rm if}\  \alpha \geq 1
                            \end{array}\right.\\
                            &=& o(1)
                          \end{eqnarray*}
                          by the choice of $k$. 
                            One also has using the assumption  (H2)
   \begin{eqnarray*}
 &|&F ((j|x_j-y_j|^{k-2} (x_j-y_j)+ {k|x_j|^{k-2} x_j\over j^{3k}},X_j + {D^2(|x|^k)\over j^{3k}}(x_j))-F ((j|x_j-y_j|^{k-2} (x_j-y_j)\\
 &+& {k|x_j|^{k-2} x_j\over j^{3k}},X_j )| \\
                             & \leq&  c j^{-6} (j |x_j-y_j|^{k-1})^\alpha  \\
                             &\leq&  \left\{ \begin{array}{cc}
                            c j^{-6-\alpha}& {\rm if}  \  \alpha <0\\
                             cj^{-6+ {\alpha\over k}} (j|x_j-y_j|^k)^{\alpha({k-1\over k})}& {\rm if} \ \alpha \geq 0
                             \end{array}\right.\\
                             &=& o(1)
\end{eqnarray*}by the choice of $k$. 
                                                         
Treating analogously  the terms involving $h$, in particular using the H\"older's regularity of $h$ with respect to $t$,  together with  (H3), one obtains
  \begin{eqnarray*}
                           &&g(x_j, t_j) \geq j(t_j-s_j) -F (j|x_j-y_j|^{k-2} (x_j-y_j)+ {k|x_j|^{k-2} x_j\over j^{3k}},X_j + {D^2(|x|^k)(x_j)\over j^{3k}})\\
                           &- & h(x_j,t_j)\cdot \left(j|x_j-y_j|^{k-2} (x_j-y_j)+ {k|x_j|^{k-2} x_j\over j^{3k }}\right) \left\vert j|x_j-y_j|^{k-2} (x_j-y_j)+ {k|x_j|^{k-2} x_j\over j^{3k}}\right\vert^\alpha  \\
                          &\geq & j(t_j-s_j)-F(x_j, j|x_j-y_j|^{k-2} (x_j-y_j),X_j )\\
                          &-&h(x_j,t_j)\cdot  j|x_j-y_j|^{k-2} (x_j-y_j)| j|x_j-y_j|^{k-1} |^\alpha - o(1)\\
                                                     &\geq &  j(t_j-s_j) -F(y_j, j|x_j-y_j|^{k-2} (x_j-y_j),-Y_j )\\
                                                     &-& h(y_j,s_j)\cdot  j|x_j-y_j|^{k-2} (x_j-y_j)| j|x_j-y_j|^{k-2} (x_j-y_j)|^\alpha |\\
                           &&-o(1) \\
                           &\geq& f(y_j, s_j) + {\kappa\over T^2}-o(1)
\end{eqnarray*}
  We now conclude as before : 
                            We use the fact that $|x_j-y_j|+ |t_j-s_j|$ goes to zero,   the uniform continuity of $f$ and $g$ ,   and we  pass to the limit to get a contradiction. 
                            
                            This ends the proof of theorem \ref{compRN}. 
                          
  \bigskip                                                          
   We now prove that the solutions are  H\"older's continuous. 
 \begin{prop}\label{prop11}
 
 Suppose that $u$ is a solution of $1_{\{f, \psi\}}$ on $\R^N\times ]0,T[$. 
 Suppose that there exist some constant   $c_1$ and $c_2$ such that 
 \begin{equation}\label{equpsi}
 u(x,t)\leq \psi (y)+ c_1|x-y|^{\gamma_\psi} + c_2 t^{1\over q(\alpha+1)-\alpha }
 \end{equation}
 
 \begin{equation}\label{equpsi2}
 u (x,t)  \geq \psi (y) -c_1|x-y|^{\gamma_\psi} - c_2 t^{1\over q(\alpha+1)-\alpha }
  \end{equation}
  We assume that $f$ is uniformly continuous and bounded, is   $\gamma_f$ H\"olderian with respect to  $t$,   uniformly in $x$,  and   that $\psi$ is H\"olderian of exponent $\gamma_\psi$ on $\R^N$  and bounded, then  $u$ is Holder's continuous  of exponent $\gamma_\psi$ with respect to $x$ and of exponent $\gamma^\star = \inf (\gamma_f, {1\over q(\alpha+1)-\alpha})$ with respect to $t$  on every compact set of $\R^N \times ]0, T[$. 
  \end{prop}
   
  We shall need the following proposition, which proves some Holder's regularity with respect to $t$, when $x$ is fixed.  
  
   \begin{prop} \label{prop12}
   
    Under the assumptions of Proposition \ref{prop11} there exists some constant $C_2$ such that for all $x\in \R^N$ and for all  $ t,s>0$
 $$|u(x,t+s)-u (x,t)|\leq C_2s^{\gamma^\star}$$
 where $\gamma^\star = \inf (\gamma_f, {1\over q(\alpha+1)-\alpha}),$ $q ={q_1\over \gamma_\psi},$ $q_1 =   {\sup (2, {\alpha+2\over \alpha+1}
)}$.

\end{prop}
{\em Proof }
   We   first use the estimates   (\ref{equpsi})  and (\ref{equpsi2}) which give for $y = x$ :
    $$|\psi (x)-u(x,s)|\leq  c_2 s^{1\over q(\alpha+1)-\alpha }$$
and   the comparison principle  in Theorem  \ref{compRN}  on $\R^N\times ]0,T[$ : 
 We define  fro $s$ fixed in $[0,T]$ and $ t\in [0, T-s]$
  $$v(x,t) = u (x, t+ s)+c_f t s^{\gamma_f}+ \sup_{x\in \R^N} |\psi(x)-u(x,s)|.$$
  where $c_f$ is some Holder's constant of $f$ with respect to $t$. 
Then $v$ is a supersolution of $1_{\{f, \psi\}}$ on $\R^N \times [0, T-s[$. 
  Let us note that  $v$ and $u$ have  the properties  
  \begin{equation}\label{nd}
  u (x,t)\leq \psi(y)+ c_1|x-y|^{\gamma_\psi}+ c_2 t^{1\over q(\alpha+1)-\alpha } \leq \psi (y)+ 2c_1(|x-y|+ 1) + c_2 T^{1\over q(\alpha+1)-\alpha }  
  \end{equation}
   and 
   \begin{equation} \label{eqbn}
   v (x,t)\geq \psi(y)- c_1|x-y|^{\gamma_\psi}- c_2 (t+s)^{1\over q(\alpha+1)-\alpha } \geq \psi (y)- 2c_1(|x-y|+ 1) -c_2 (2T)^{1\over q(\alpha+1)-\alpha }  
   \end{equation}
   and $u(x,0)\leq v(x,0)$ by construction.

   Hence one can apply the comparison theorem \ref{compRN} to obtain that 
  $$u (x, t) \leq v (x,t)+ \sup_{x\in \Omega} |u  (x,0)-v (x,0)|\leq u (x, t+s) + c_f T s^{\gamma_f}+ c_2 s^{1\over q(\alpha+1)-\alpha }  $$
                                     In the same manner defining 
                                   $v (x, t)= u (x, t+s)-c_f t s^{\gamma_f}-\sup_x |\psi(x)-u(x,s)|$
                                   then $u (x, t)$ and $v$ are super and sub-solution for the same equation, and then using theorem \ref{compRN} one gets 
                                   $$u (x, t)\geq  u (x,t+s)-c_f T s^{\gamma_f}-c_2 s^{1\over q(\alpha+1)-1}.$$
                                    The result follows.

                 \bigskip
     Proof of  proposition \ref{prop11}  
     
     First  we observe that $u$ is bounded, taking $y= x$ in the inequalities (\ref{equpsi}) and (\ref{equpsi2})
       and using the fact that $\psi$ is bounded. 
                                 
 Let $\delta$ be given less than $1$,     $L > \sup (4c_1 + L_\psi ,  \left({4 |f|_\infty\over \gamma_\psi ^{1+\alpha} ( 1-\gamma_\psi) }\right)^{1\over 1+\alpha})$ and $ M \geq \sup ({2\sup u\over  \delta^{\gamma^\star}}, c_2, {2 c_2 T^{\gamma^\star}\over \delta^{\gamma^\star}})$.    We define the set
  $$\Delta_\delta = \{ (x,y, t,s), |x-y|< \delta ,  |t-s|< \delta, (t,s)\in ]0,T[\}$$
 and for $j$ large the function 
$$ \psi_j (x,y, t,s) = u(x,t)-u(y,s)-L|x-y|^{\gamma_\psi} -{|x|^2\over 2j^2} -M  |t-s|^{\gamma^\star}.$$
 We shall prove that for $j$ large enough, $\psi_j$ is $\leq 0$. The result   will follow 
 by  passing to the limit on   each compact set of $\R^N\times ]0,T[$.                                      

  We then assume  by contradiction that $\psi_j$ has a maximum strictly positive. 
   Then for  $\kappa$ small enough 
   $$\psi_j-{\kappa\over T-t} -{\kappa\over T-s}$$
   has also its supremum strictly positive and 
       we begin  to observe that on the boundary of $\Delta_\delta$, this function is $\leq 0$. 
 
  Indeed in the case where  $|t-s| = \delta$  then  by hypothesis (\ref{equpsi}) and (\ref{equpsi2})
 
$$u(x,t)-u(y,s) \leq c_1 |x-y|^{\gamma_\psi} + {2 c_2 T^{\gamma^\star}}\leq L|x-y|^{\gamma_\psi} + M |t-s|^{\gamma^\star}$$
 In the case where $t = 0$, $s>0$ and $|x-y|\leq \delta $ one uses once   more (\ref{equpsi}) and (\ref{equpsi2}). 
 
  Finally    the  supremum cannot be achieved for $t= T$ or $s = T$.

                                   Let us note that if $\psi_j$ has a supremum $>0$, 
                                   $$\psi_j ^n(x,t,y,s) = u(x,t)-u(y,s)-L |x-y|^{\gamma_\psi} -{|x|^2\over 2j^2} -{M} ({1\over n^2} + [t-s|^2)^{\gamma^\star\over 2}-{\kappa\over T-t} -{\kappa\over T-s}$$
                                    has also a supremum $>0$  achieved inside $\Delta_\delta $, for $n$ large enough.  We fix $n$ large enough. Let  
                                    $(x_j, y_j, t_j, s_j)$  be a point where the supremum of $\psi_n$ is achieved.   By the previous considerations, it cannot be achieved on the boundary.  By proposition \ref{prop12} one has $x_j\neq y_j$ and then the function $|x-y|^{\gamma_\psi}$ is ${\cal C}^2$ on a neighborhood of $(x_j, y_j)$. 
Using Ishii's lemma (see also Lemma 2.1 in \cite{BD1} ) we have the existence of $(X_j, Y_j)$ with 
                                   \begin{eqnarray*}
&&\left({ \gamma^\star  M(t_j-s_j) ({1\over n^2} + |t_j-s_j|^2)^{1-{\gamma^\star\over 2}}}
                                                     + {\kappa\over (T-t_j)^2},\right. \\
                                                      &&\left.\gamma_\psi L(x_j-y_j)|x_j-y_j|^{\gamma_\psi-2} +  {x_j\over j^2},
                                                         X_j+ {I\over  j^2} \right)\in J^{2,+} u(x_j, t_j)
                                 \end{eqnarray*}

    \begin{eqnarray*}      ({\gamma^\star M(t_j-s_j) ({1\over n^2} + |t_j-s_j|^2)^{1-{\gamma^\star\over 2}}}&-&{\kappa\over (T-s_j)^2}, \gamma_\psi L(x_j-y_j)|x_j-y_j|^{\gamma_\psi-2} ,  -Y_j)\\
     &\in& J^{2,-} u(y_j, s_j)
      \end{eqnarray*}

          with 
      $$\left(\begin{array}{cc}
      X_j&0\\
      0& Y_j
      \end{array}\right)\leq \left(\begin{array} {cc} 
      B(x_j, y_j)&-B(x_j, y_j)\\
      -B(x_j, y_j)&B(x_j, y_j)\end{array}\right)$$
      
      with $B(x,y) =L \gamma_\psi |x-y|^{\gamma_\psi-2} (I+ (\gamma_\psi-2){(x-y)\otimes (x-y)\over |x-y|^2}) = D^2 (|X|^{\gamma_\psi})(x-y)$
      
      Let us observe that  due to the hypothesis,  
      $|{x_j\over j^2}|\leq {c\over j} \leq  {\gamma_\psi \over 2 L  \delta ^{\gamma_\psi-1} }$, 
   and then 
       $|\gamma_\psi L(x_j-y_j)|x_j-y_j|^{\gamma_\psi-2} +  {x_j\over j^2}|\geq { {\gamma_\psi\over 2} L |x_j-y_j|^{\gamma_\psi-1}}$.

         We use as in the proof of theorem 2, the inequality          $$|tr(X_j+ Y_j)|= -tr(X_j+Y_j) \geq  {4\gamma_\psi} (1-\gamma_\psi)L |x_j-y_j|^{\gamma_\psi-2} $$
         and  the  fact that  for some constant $c$

         $$|X_j|+ |Y_j|\leq c (|tr(X_j+ Y_j)|$$

          We then  use the property (H6)  of $F$ to get that 
        \begin{eqnarray*}
        |Fx_j,    \gamma_\psi L(x_j-y_j)|x_j-y_j|^{\gamma_\psi-2} &&+  {x_j\over j^{2}},  X_j+ {I\over  j^{2}} )\\
               &-&F(x_j,  \gamma_\psi L (x_j-y_j)|x_j-y_j|^{\gamma_\psi-2} , X_j) |\\
        &\leq&  O(j^{-1}) ( L|x_j-y_j|^{\gamma_\psi-1})^{\alpha-1} |X_j|+ O({1\over j^2})( L|x_j-y_j|^{\gamma_\psi-1})^{\alpha}   \\
        &\leq &o(1 ) (L|x_j-y_j|^{\gamma_\psi-1})^{\alpha}| tr(X_j+ Y_j)|. 
        \end{eqnarray*}
        
        And we use only the fact that $h$ is bounded to observe  that 
        $$|h(x_j, t_j)-h(y_j, s_j)\cdot  (\gamma_\psi L)^{1+\alpha}(x_j-y_j)|x_j-y_j|^{(1+\alpha)(\gamma_\psi-1)-1}|\leq o(1)(\gamma_\psi  L |x_j-y_j|^{\gamma_\psi-1})^\alpha |tr(X_j+Y_j)|$$

       We now write 
       \begin{eqnarray*}
        f(x_j, t_j)&\geq &  {\gamma^\star  M(t_j-s_j)({1\over n^2} + |t_j-s_j|^2)^{1-{\gamma^\star\over  2}}}+ {\kappa \over( T-t_j)^2} 
        \\
        &-&F(x_j, \gamma_\psi L(x_j-y_j)|x_j-y_j|^{\gamma_\psi-2} + {x_j\over j^{2}},  X_j+ {I\over  j^{2}} )\\
        &-&h(x_j, t_j)\cdot  (\gamma_\psi L)^{1+\alpha}(x_j-y_j)|x_j-y_j|^{(1+\alpha)(\gamma_\psi-1)-1}\\
        &\geq &  {\gamma^\star M(t_j-s_j)({1\over n^2} + |t_j-s_j|^2)^{1-{\gamma^\star\over 2}}} -{\kappa\over (T-s_j)^2}\\
        &
        -&F(y_j,  \gamma_\psi L(x_j-y_j)|x_j-y_j|^{\gamma_\psi-2} ,  -Y_j)\\
        &-&h(y_j, s_j)\cdot  (\gamma_\psi L)^{1+\alpha}(x_j-y_j)|x_j-y_j|^{(1+\alpha)(\gamma_\psi-1)-1} \\
         &&+ (\gamma_\psi L|x_j-y_j|^{\gamma_\psi-1})^{\alpha} tr(X_j+Y_j) + o(1) |\gamma_\psi L |x_j-y_j|^{\gamma_\psi-1}|^\alpha (|tr(X_j+Y_j)|)\\
        &\geq & f(y_j, s_j)+  (\gamma_\psi L|x_j-y_j|^{\gamma_\psi-1})^{\alpha} tr(X_j+Y_j) (1- o(1)).\\
        \end{eqnarray*}
        
        We have obtained a contradiction since  this would imply  that 
        $$      (  \gamma_\psi L|x_j-y_j|^{{\gamma_\psi}-1})^{\alpha} L|x_j-y_j|^{
      \gamma_\psi-2} (1-o(1))\leq 2|f|_\infty, $$
        
      which is absurd  by  the choice of the constant $L$.

      This ends the proof of the following Holder's result :
         
          \begin{prop}
          
          Suppose that $\psi$ is H\"olderian and bounded in  $\R^N$ and that  $f$ is uniformly continuous and Holder's in $t$,  uniformly w.r.t. $x$.  Then, 
          there exists a  unique viscosity solution of $1_{\{f, \psi\}}$ on $\R^N\times ]0,T[$. This solution is Holder's continuous on every compact set of $\R^N\times ]0,T[$. 
          \end{prop}
          Hence using Ascoli's theorem,  we have also  
          \begin{cor}
           Let $(f_n$, $\psi_n)$  be a sequence of bounded Holder's continuous functions,  $\psi_n$  being uniformly Holder's and $(f_n)$ being  uniformly Holder's in $t$,  uniformly w.r.t. $x$.  Then the sequence $(u_n)$  of solutions of $1_{\{f_n, \psi_n\}}$ is relatively compact on  every compact set of $ \R^N\times ]0,T[$.            \end{cor}

       \section{Appendix}
       
       In this appendix we prove that the solutions of Ohnuma and Sato in the case where $\alpha \neq 0$ are the same as our solutions. In the same manner we prove that it is also the case for  the infinity Laplacian  using the adapted definition  of Evans and Spruck, and Juutinen and Kawhol.

              \subsection{The case $\alpha \neq 0$}
       
     The reader can consult \cite{OS}  for the definition of ${\cal F} (F)$ and ${\cal A}(F)$. 
     
     We recall that in \cite{OS} the right hand side  $f$ is zero. 
           \begin{prop}
           
            The solutions in our sense are the same as the solutions in the Ohnuma and Sato's sense. 
            \end{prop}
            Proof 
           
            Suppose that $u$ is a supersolution of $1_{\{0\}}$ in the Ohnuma and Sato's sense. 
             Suppose that $(\bar x, \bar t)$ is some point such that   for some $\delta_1$  and for some ${\cal C}^1$ function $\varphi$ on $]0,T[$ :    
             $$\inf_{|t-\bar t|< \delta_1}  (u(\bar x, t)-\varphi(t)) =   u(\bar x, \bar t)-\varphi(\bar t)=0$$
             and such that 
             $x\mapsto \inf_{|t-\bar t|< \delta_1}( u(x, t)-\varphi(t))$ is 
             constant on $B(\bar x, \delta)$ for some $\delta>0$ . Then
             in particular 
             $$\inf_{x\in B(\bar x , \delta), |t-\bar t|< \delta_1}  (u(x, t)-\varphi(t))$$ has its infimum  equals to zero achieved on $(\bar x, \bar t)$. Then,   for $\epsilon >0$ the function 
             $$h(x,t) = \varphi(\bar t) + \varphi^\prime (\bar t) (t-\bar t)+{1\over 2} (\varphi^"(\bar t)-\epsilon)  (t-\bar t)^2$$
             which  belongs to ${A(f)} $,   \cite{OS}, satisfies   
             $$ \inf_{(|t-\bar t|< \delta_1, x\in B(\bar x, \delta)} (u (x,t)-h(x,t)) =0$$ Indeed 
             $$\inf_{|t-\bar t|< \delta_1, x\in B(\bar x, \delta)}( u-h)(x,t) \leq u(\bar x, \bar t)-\varphi(\bar t)=0.$$
             Moreover   for $t$ close to $\bar t$
             $$\varphi (t) \geq \varphi(\bar t) + \varphi^\prime (\bar t) (t-\bar t) + {1\over 2} (\varphi^"(\bar t)-\epsilon)  (t-\bar t)^2$$
             hence 
             $$\inf_{(|t-\bar t|< \delta_1, x\in B(\bar x, \delta)}(u-h)(x,t)\geq\inf_{|t-\bar t|< \delta_1, x\in B(\bar x, \delta)} (u(x,t)-\varphi(t))$$
             and then   since $u$ is a supersolution of $1_{\{0\}}$, 
             $\varphi^\prime (\bar t) \geq 0$ which is the desired conclusion. 
             
             We want to prove the reverse sense. We assume that $u$ is a super solution in our sense. We suppose that $(\bar x, \bar t)$ and $\varphi$ are such that 
             $(u-\varphi) \geq (u-\varphi )(\bar x, \bar t )= 0$, 
             with $\varphi \in {\cal A} (F)$. 
             
           Let $f\in {\cal F} (F)$ and $\omega$ be  a continuous function such that $\omega (0)=0, $ $\omega (t-\bar t) = o(|t-\bar t|)$, be such that for $(x,t) \in V$ a neighborhood of $(\bar x, \bar t)$,  
           $$|\varphi (x, t)-\varphi (\bar x, \bar t)-\partial_t \varphi (\bar x, \bar t)(t-\bar t)|\leq f(|x-\bar x|) + \omega (t-\bar t)$$
           Then 
          $$ h(x,t) : = \varphi (\bar x, \bar t)+\partial_t \varphi (\bar x, \bar t)(t-\bar t)-  f(|x-\bar x|) - \omega (t-\bar t)\leq \varphi (x,t)$$
          Moreover 
           
            $$\inf_{(x,t)\in V} (u(t,x) - h(x,t)) = 0$$
             Indeed 
            $$\inf _{(x,t)\in V}(u(x,t)-h(x,t)) \leq u(\bar x, \bar t)-h(\bar x, \bar t)$$
            secondly  by the previous remark, 
            $$u-h\geq u-\varphi .$$
            
             Now acting as in lemma \ref{lem1}   ie replacing $C_1|x-\bar x|^k$ by $f(|x-\bar x|)$ and $C_2 |t-\bar t|^2$ by
 $\omega (|t-\bar t|)$
  one gets since $\lim_{x\rightarrow 0}F(\nabla  f, D^2 f)(|x|)=0$   that 
              $\partial_t \varphi (\bar x , \bar t)\geq 0$, 
              which is the desired conclusion. 
              
              \subsection{ The case $\alpha= 0$ and the infinity Laplacian}
              
               We  prove here that  our definition is equivalent to the one of Evans and Spruck in the case of the infinity Laplacian (see also \cite{JK}).

          We shall need  the following lemma, whose proof is postponed  for the sake of clearness.                 
 \begin{lemme}\label{lem2} 
 
 Suppose that $u$ is a supersolution of 
 $$u_t-F(x, \nabla u, D^2 u)-h(x,t)\cdot \nabla u|\nabla u|^\alpha \geq f(x,t)$$
 and suppose that  $\varphi$ is some ${\cal C}^2$ function on $]0,T[$,  with $\varphi(\bar t)=0$,  that $k > \sup (2, {\alpha+2\over \alpha+1})$,   that $M$ is some symmetric matrix 
 and  $(0, \bar t)\in \Omega \times ]0, T[$ are such that  for some  $\delta_1>0$
 $$\inf _{x\in B(0, \delta_1), | t-\bar t|< \delta_1}  (u(x,t)-\varphi(t)-{1\over 2}( Mx, x) )=  u(0, \bar t)$$
 Then 
 $$\varphi^\prime (\bar t) -{\cal M}_{a,A}^-(M)\geq f(0, \bar t).$$
  
 \end{lemme}
 
 We postpone the proof of Lemma \ref{lem2}
 
 We now consider a supersolution $u$ in our sense and assume that $\varphi$ is some ${\cal C}^2$ function which achieves $u$ by below on $(\bar x, \bar t)$ with $\nabla _x \varphi (\bar x, \bar t) = 0$.  
We apply  lemma \ref{lem2}  with  $\bar x$ in place of $0$,  $\nabla_x \varphi(\bar x,\bar  t)= 0$ and replacing $\varphi(t)$by $ \partial_t \varphi (\bar x,\bar t) (t-\bar t) $,  and  $M = {D^2\varphi} (\bar x, \bar t)$ 
  one gets the desired conclusion. 
  
  Proof of lemma \ref{lem2}:

  For $C_2>0$  one still has 
   $$\inf _{x\in B(0, \delta_1), | t-\bar t|< \delta_1}  (u(x,t)-\varphi(t)-{1\over 2}( Mx, x) + C_2(t-\bar t)^2)=  u(0, \bar t)$$
  and the infimum is strict in $t$. 
  
   We assume  first that 
  $ x\mapsto \inf _{|t-\bar t|<  \delta_1}  (u(x,t)-\varphi(t)+ C_2(t-\bar t)^2 )$ is equal to $u(0, \bar t)$ and is constant w.r.t.  $x$ in a neighboorhood of  $\bar x$. We then prove that 
  $M \leq  0$ and 
  $\varphi^\prime (\bar t)\geq f(0, \bar t)$.
  
  This will imply that 
  $\varphi^\prime (\bar t) -{\cal M}_{a,A}^-(M) \geq f(0, \bar t)$.
  
  Indeed one has  for all $x$ in  a neighborhhod of $0$, 
  $u(0, \bar t) = \inf _{|t-\bar t|< \delta_1}  (u(x,t)-\varphi(t)+ C_2(t-\bar t)^2) $
  and  also by hypothesis 
  $$u(0, \bar t) = \inf _{(|t-\bar t|<  \delta_1), x\in B(0, \delta_1)}  \{u(x,t)-\varphi(t)-{1\over 2}(Mx,x) +C_2(t-\bar t)^2\}$$
   and   then for all $x$ in  a neighborhhod of $0$, 
 $$u(0, \bar t) \leq  \inf _{|t-\bar  t|<  \delta_1}  \{u(x,t)- \varphi(t)+ C_2(t-\bar t)^2 \}-{1\over 2}(Mx,x)= u(0, \bar t)-{1\over 2}(Mx,x)$$
    This implies that 
     for all $x$ in  a neighborhhod of $0$, 

     $$(M x,x)\leq 0, $$ or equivalently that $M$ is  a nonpositive symmetric matrix. 
     Using the definition,  as we pointed out before, 
       $\varphi^\prime (\bar t)\geq f(0, \bar t)$ and this  implies  the desired result. 
       
       We now assume that we are not in the case  where 
       $ x\mapsto \inf _{|t-\bar t|<  \delta_1)}  (u(x,t)-\varphi(t)+ C_2(t-\bar t)^2 )$ is equal to $u(0, \bar t)$ and is constant w.r.t.  $x$ in a neighboorhood of  $\bar x$.
       
  For the sequel one can assume that $M$ is invertible. indeed, if it is not the case there exists $\epsilon>0$ arbitrarily small in order that $M-\epsilon Id$ is invertible.  Moreover $M-\epsilon Id$ is also such that   
    $$ \inf _{(|t-\bar t|<  \delta_1), x\in B(0, \delta_1)}  \{u(x,t)-\varphi(t)-{1\over 2}((M-\epsilon Id)(x),x) +C_1|x|^k +C_2(t-\bar t)^2\} = u(0, \bar t)$$
    So we shall prove that 
    $$\varphi^\prime (\bar t)-{\cal M}_{a,A}^-( M-\epsilon Id) \leq  f(0, \bar t)$$ and we shall get the result by passing to the limit with $\epsilon$. 
    
     So from now we assume that $M$ is invertible. 
      
     For $k>2 $ and for  all positive constant $C_1$   then 
     $$ \inf _{(|t-\bar t|<  \delta_1), x\in B(0, \delta_1)}  \{u(x,t)-\varphi(t)-{1\over 2}(Mx,x) +C_1|x|^k +C_2(t-\bar t)^2\}$$
      has also its  infimum achieved on $(0, \bar t)$,  and this  infimum is strict in $x$ and $t$.  Hence for all $\delta >0$  there exists $\epsilon (\delta)>0$ such that 
  \begin{eqnarray*}
  \inf\left(  \inf _{(|t-\bar t|> \delta, x\in B(0, \delta_1)}  \right.&&\{u(x,t)-\varphi(t)-{1\over 2}(Mx,x) 
  +C_1|x|^k +C_2(t-\bar t)^2\} 
  \\
  &&\left. \inf _{(|t-\bar t|<  \delta_1, |x|> \delta }  \{u(x,t)-\varphi(t)-{1\over 2}(Mx,x) +C_1|x|^k +C_2(t-\bar t)^2\}\right)\\
&    >& u(0, \bar t)+ \epsilon (\delta)
\end{eqnarray*}
   
In the following we choose $\delta$ such that $(2\delta)^{k-1} < {\inf_{\lambda_i \in Sp (M)}|\lambda_i(M)| \over 2kC_1}$.  
 Let then $\delta_2 $  be such that $\delta_2 < \delta $ and 
 $$k(2\delta_1)^{k-1} C_1 \delta_2 + |M|_\infty (\delta_2^2 + 2\delta_2\delta_1 ) \leq {\epsilon/4}$$
 
 With this choice,   using the fundamental calculus theorem, one gets that  for $x\in B(0,   \delta_2)$,  
 
   \begin{eqnarray}\label{eq1}
    \inf _{\{|t-\bar t|< \delta), y\in B(0, \delta)\}} && \{u(y,t)-\varphi(t)-{1\over 2}(M(y-x),(y-x)) +C_1|x-y|^k +C_2(t-\bar t)^2\} \nonumber \\
    &\leq&   \inf _{(|t-\bar t|< \delta_1), y\in B(0, \delta_1)}  \{u(y,t)-\varphi(t)-{1\over 2}(My,y) +C_1|y|^k +C_2(t-\bar t)^2\}+ {\epsilon \over 4}\nonumber \\
    &=& u(0, \bar t)+ {\epsilon\over 4}
    \end{eqnarray}
     while

  \begin{eqnarray}\label{eq2}
  \inf \left( \inf _{\{|t-\bar t|< \delta_1), |y|> \delta \}}  \right.&&(\{u(y,t)-\varphi(t)-{1\over 2}(M(y-x),(y-x)) +C_1|x-y|^k +C_2(t-\bar t)^2\}, \nonumber\\
  &&\left.  \inf _{(|t-\bar t|> \delta, y\in B(0, \delta_1)}  \{u(y,t)-\varphi(t)-{1\over 2}(M(y-x),(y-x))\right.\nonumber \\
  &+&\left.C_1|x-y|^k +C_2(t-\bar t)^2\} \right)\nonumber\\
  &\geq &u(0, \bar t)+ {3\epsilon\over 4}
  \end{eqnarray}

  We choose $x_\delta$ as follows  : 
  Since the function $ \inf _{|t-\bar  t|<  \delta_1)}  (u(x,t)-\varphi(t) +  C_2 |t-\bar t|^2) $ is not constant around $\bar x$, for all $\delta >0$ there exists $x_\delta$ and $y_\delta$  in $B(0, \delta_2) $ such that 
  \begin{eqnarray*}
   \inf_{|t-\bar t|< \delta_1}\{u(x_\delta ,t)-\varphi(t) &+&  C_2 |t-\bar t|^2\} \\
   &> & \inf_{|t-\bar t|< \delta_1}\{u(y_\delta ,t)-\varphi(t) +  C_2 |t-\bar t|^2-\left({1\over 2}(M(x_\delta-y_\delta), x_\delta -y_\delta) \right)\\
   &+& C_1 |x_\delta -y_\delta|^k\}
   \end{eqnarray*}
  
  Then the infimum 
  $ \inf _{(|t-\bar t|< \delta_1), y\in B(0, \delta_1)}  (u(y,t)-\varphi(t)-{1\over 2}(M(y-x_\delta ),(y-x_\delta )) +C_1|x_\delta -y|^k +C_2(t-\bar t)^2),$ is  achieved on some point $(z_\delta, t_\delta)$ with $z_\delta \neq x_\delta$. Moreover by ( \ref{eq1}) and (\ref{eq2}) the infimum is achieved in $B(0, \delta)\times ]\bar t-\delta, \bar t+\delta[$. Let $(z_\delta, t_\delta)$ be a point on which this infimum is achieved,  
  then 
  \begin{eqnarray*}
  \psi (x,t)& =& \varphi(t)+ {1\over 2}(M(x-x_\delta), x-x_\delta) )-  {1\over 2}(M(z_\delta-x_\delta),(z_\delta-x_\delta)) \\
  &+& C_1 |x_\delta -z_\delta|^k-C_1 |x_\delta -x|^k \\
  &-&C_2(t-\bar t)^2+ C_2 (t_\delta -\bar t)^2
  \end{eqnarray*}
  
achieves 
  $u$  by below on $(z_\delta, t_\delta)$.  
  
  With the choice of $\delta$,  the gradient of $\psi$  on $z_\delta$, which equals 
  $M(z_\delta-x_\delta) + kC_1 |x_\delta-z_\delta|^{k-2} (x_\delta -z_\delta) $ is different from zero, since  $z_\delta \neq x_\delta$.  Indeed if it was the case, $x_\delta-z_\delta$ would be an eigenvector for $M$ corresponding to the eigenvalue $kC_1 |x_\delta-z_\delta |^{k-1}$, which is impossible since 
  $kC_1 (2\delta)^{k-1} < \inf_i (|\lambda_i (M)|$. 
  Using the fact that $u$ is a supersolution one gets  that 
   \begin{eqnarray*}
   \varphi^\prime (t_\delta) &-& F(M(z_\delta-x_\delta) + kC_1 |x_\delta-z_\delta|^{k-2} (x_\delta -z_\delta) , M 
  -C_1 D^2 (|x_\delta-z|^k)(z_\delta)) \\
  &-&h(z_\delta , t_\delta) \cdot \nabla \psi (z_\delta, t_\delta)| \nabla \psi (z_\delta, t_\delta)^\alpha \\
  &\geq& f(z_\delta, t_\delta)
  \end{eqnarray*}
  
      and then 
      \begin{eqnarray*}
      \varphi^\prime (t_\delta) &-& {\cal M}_{a,A} ^-(M -C_1 D^2 (|x_\delta-z|^k)(z_\delta)\\
&       -&h(z_\delta , t_\delta) \cdot  \nabla \psi (z_\delta, t_\delta)| \nabla \psi (z_\delta, t_\delta)|^\alpha \\
&\geq& f(z_\delta, t_\delta)
\end{eqnarray*}

   Letting $\delta$  go to zero and using $z_\delta \in B(0, \delta_2)\subset B(0, \delta)$, $|t-t_\delta|< \delta $,   $k>2 $ and the lower semicontinuity of $f$ one gets 
   $$\varphi^\prime (\bar t) -{\cal M}^-_{a, A} (M) \geq f(0, \bar t)$$

 {\em Acknowledgment : 
 The author is very grateful to Isabeau Birindelli for the interest she brought to this  subject and the precious advices she gave, which permit to improve it.}

\end{document}